\documentclass[12pt]{article}

\usepackage{amsmath,amssymb}

\newcommand{\non}{\nonumber}
\newcommand{\scl}{\scriptstyle}
\newcommand{\wt}{\widetilde}
\newcommand{\wh}{\widehat}
\newcommand{\ot}{\otimes}

\newcommand{\ts}{\,}
\newcommand{\tss}{\hspace{1pt}}
\newcommand{\U}{ {\rm U}}
\newcommand{\Z}{ {\rm Z}}
\newcommand{\Sr}{ {\rm S}}
\newcommand{\Y}{ {\rm Y}}
\newcommand{\SY}{ {\rm SY}}
\newcommand{\Pf}{ {\rm Pf}\ts}
\newcommand{\Hf}{ {\rm Hf}\ts}
\newcommand{\tr}{ {\rm tr}\ts}
\newcommand{\CC}{\mathbb{C}}

\newcommand{\Ar}{ {\rm A}}
\newcommand{\Ir}{ {\rm I}}
\newcommand{\End}{{\rm{End}\ts}}
\newcommand{\Sym}{\mathfrak S}
\newcommand{\sgn}{ {\rm sgn}\ts}
\newcommand{\tpr}{t^{\tss\prime}}
\newcommand{\spr}{s^{\tss\prime}}
\newcommand{\h}{\mathfrak h}
\newcommand{\gl}{\mathfrak{gl}}
\newcommand{\oa}{\mathfrak{o}}
\newcommand{\spa}{\mathfrak{sp}}
\newcommand{\g}{\mathfrak{g}}
\newcommand{\agot}{\mathfrak{a}}
\newcommand{\sll}{\mathfrak{sl}}
\newcommand{\la}{\lambda}
\newcommand{\La}{\Lambda}

\newcommand{\om}{\omega^{}}
\newcommand{\qdet}{ {\rm qdet}\ts}
\newcommand{\sdet}{ {\rm sdet}\ts}
\newcommand{\gr}{ {\rm gr}}
\newcommand{\Proof}{\noindent{\it Proof.}\ \ } 
\newcommand{\Outline}{\noindent{\it Outline of the proof.}\ \ }

\newtheorem{thm}{Theorem}[section]
\newtheorem{prop}[thm]{Proposition}
\newtheorem{lem}[thm]{Lemma}
\newtheorem{cor}[thm]{Corollary}
\newtheorem{defin}[thm]{Definition}
\newtheorem{example}[thm]{Example}

\newcommand{\bth}{\begin{thm}}
\renewcommand{\eth}{\end{thm}}
\newcommand{\bpr}{\begin{prop}}
\newcommand{\epr}{\end{prop}}
\newcommand{\ble}{\begin{lem}}
\newcommand{\ele}{\end{lem}}
\newcommand{\bco}{\begin{cor}}
\newcommand{\eco}{\end{cor}}
\newcommand{\bex}{\begin{example}}
\newcommand{\eex}{\end{example}}
\newcommand{\bde}{\begin{defin}}
\newcommand{\ede}{\end{defin}}
\newcommand{\bal}{\begin{aligned}}
\newcommand{\eal}{\end{aligned}}
\newcommand{\beq}{\begin{equation}}
\newcommand{\ben}{\begin{equation*}}

\def\beql#1{\begin{equation}\label{#1}}

\newbox\squ  
\setbox\squ=\hbox{\vrule width.3pt
            \vbox{\hrule height.3pt width.4em\kern1ex\hrule height.3pt}%
            \vrule width.3pt}
\def\endproof{%
  \ifmmode\eqno\copy\squ\smallskip\else{\unskip\nobreak\hfil%
    \penalty50\hskip2em\hbox{}\nobreak\hfil\copy\squ
    \parfillskip=0pt \finalhyphendemerits=0\penalty-100\smallskip}
  \fi}

\begin{document}

\title{\Large\bf 
Yangians and their applications\thanks{\noindent
To 
appear in {\it Handbook of Algebra\/} (M. Hazewinkel, Ed.), 
Vol. 3, Elsevier.
}}
\author{{\sc A. I. Molev}\\[15pt]
School of Mathematics and Statistics\\
University of Sydney,
NSW 2006, Australia\\
{\tt
alexm\hspace{0.09em}@\hspace{0.1em}maths.usyd.edu.au
}
}

\date{}
\maketitle


\newpage

\tableofcontents

\newpage

\pagestyle{plain}
\setcounter{page}{1}

\section{Introduction}\label{sec:int}
\setcounter{equation}{0}

The term {\it Yangian\/} was introduced by Drinfeld (in honor of C.~N.~Yang)
in his fundamental paper \cite{d:ha} (1985).
In that paper, besides the Yangians, Drinfeld defined the {\it quantized
Kac--Moody algebras\/} which together with the work
of Jimbo~\cite{j:qd}, who introduced these algebras independently,
marked the beginning of the era of {\it quantum groups\/}.
The Yangians form a remarkable family of
quantum groups related to rational solutions of the classical
Yang--Baxter equation. For each simple finite-dimensional Lie
algebra $\agot$ over the field of complex numbers, the corresponding
Yangian $\Y(\agot)$ is defined as a canonical deformation 
of the universal enveloping algebra $\U(\agot[x])$ for the polynomial current
Lie algebra $\agot[x]$. Importantly, the deformation is considered
in the class of Hopf algebras which guarantees its uniqueness under some
natural `homogeneity' conditions.
An alternative description of the algebra $\Y(\agot)$ was given
later in Drinfeld~\cite{d:nr}.

Prior to the introduction of the Hopf algebra $\Y(\agot)$ in \cite{d:ha}, the
algebra, which is now called the {\it Yangian for the general 
linear Lie algebra\/}
$\gl_n$ and denoted $\Y(\gl_n)$, was considered in the work of Faddeev
and the St.-Petersburg school in relation with the {\it inverse
scattering method}; see for instance
Takhtajan--Faddeev~\cite{tf:qi}, Kulish--Sklyanin~\cite{ks:qs},
Tarasov~\cite{t:sq, t:im}.
The latter algebra is a deformation of the
universal enveloping algebra $\U(\gl_n[x])$.

For any simple Lie algebra $\agot$ the Yangian $\Y(\agot)$ 
contains the universal enveloping
algebra $\U(\agot)$ as a subalgebra. However, only in the case
$\agot=\sll_n$ does there exist an evaluation homomorphism
$\Y(\agot)\to\U(\agot)$ identical on the subalgebra $\U(\agot)$; see 
Drinfeld~\cite[Theorem~9]{d:ha}.
In this chapter we concentrate on this distinguished
Yangian which is closely related to $\Y(\gl_n)$. 
For each of the classical Lie algebras
$\agot=\oa_{2n+1},\ts\spa_{2n},\ts\oa_{2n}$ Olshanski \cite{o:ty} introduced
another algebra	called the {\it twisted Yangian\/}.
Namely, the Lie algebra $\agot$ can be regarded as a fixed point subalgebra
of an involution $\sigma$
of the appropriate general linear Lie algebra $\gl_N$.
Then the twisted Yangian $\Y(\gl_N,\sigma)$ can be defined as a subalgebra
of $\Y(\gl_N)$
which is a deformation of the universal enveloping
algebra for the twisted polynomial current Lie algebra
\beql{twpolc}
\gl_N[x]^{\sigma}=\{\ts
A(x)\in\gl_N[x]\ts|\ts\sigma\bigl(A(x)\bigr)=A(-x)\ts\}.
\end{equation}
For each $\agot$ the twisted Yangian contains $\U(\agot)$ as a subalgebra,
and an analog of the evaluation homomorphism $\Y(\gl_N,\sigma)\to\U(\agot)$
does exist.	Moreover, the twisted Yangian turns out to be a (left) coideal
of the Hopf algebra $\Y(\gl_N)$.

The defining relations of the Yangian $\Y(\gl_n)$ can be written
in a form of a single {\it ternary\/} (or $RTT$) relation
on the matrix of generators. 
This relation has a rich and extensive background. It originates from
the quantum inverse scattering theory;
see e.g. Takhtajan--Faddeev~\cite{tf:qi}, Kulish--Sklyanin~\cite{ks:qs},
Drinfeld~\cite{d:qg}. The Yangians were primarily regarded as a vehicle for
producing rational solutions of the Yang--Baxter equation which plays
a key role in the theory integrable models;
cf. Drinfeld~\cite{d:ha}. Conversely, the ternary
relation was used in Reshetikhin--Takhtajan--Faddeev~\cite{rtf:ql} as a tool
for studying quantum groups. Moreover, the Hopf algebra structure
of the Yangian can also be conveniently described in a matrix form. 

Similarly, the twisted Yangians can be equivalently presented by generators
and defining relations which can be written
as a {\it quaternary\/} (or {\it reflection\/}) equation
for the matrix of generators, together with a {\it symmetry\/} relation.
Relations of this type appeared for the first time in Cherednik~\cite{c:fp} and
Sklyanin~\cite{s:bc}, where integrable systems with boundary conditions were
studied.

This remarkable form of the defining relations for the Yangian and twisted Yangians
allows special algebraic techniques (the so-called $R$-{\it matrix formalism})
to be used to describe the algebraic structure and study 
representations of these algebras. On the other hand, the evaluation homomorphisms
to the corresponding classical enveloping algebras allow one to use these results 
to better understand the classical Lie algebras themselves. In particular, new
constructions of the Casimir elements can be obtained in this way.
These include the noncommutative characteristic polynomials
for the generator matrices and the Capelli-type determinants.
Some other applications include the constructions of Gelfand--Tsetlin bases
and commutative subalgebras. Moreover, as was shown by Olshanski~\cite{o:ri, o:ty},
the Yangian and the twisted Yangians can be realized as some projective limit
subalgebras of a sequence of centralizers in the classical enveloping algebras.
This is known as the {\it centralizer construction\/}; see also \cite{mo:cc}.

The representation theory of the Yangians $\Y(\agot)$ is a very much
nontrivial and fascinating area. Although the 
finite-dimensional irreducible representations of $\Y(\agot)$
are completely described by Drinfeld~\cite{d:nr}, their general structure
still remains unknown. This part of the theory of the Yangians
will have to be left outside this chapter. We give, however, some references
in the bibliography which we hope cover at least some of the most 
important results in the area.

The Yangians, as well as their super and $q$-analogs,
have found a wide variety of applications in physics.
This includes the theory of integrable models in statistical mechanics,
conformal field theory,	quantum gravity. We do not attempt to review
all the relevant physics literature but some references
are given as a guide to such applications.

At the end of each section we give brief bibliographical comments
pointing towards the original articles and to
the references where the proofs or further details
can be found.
The author would like to thank
A.~N.~Kirillov, 
M.~L.~Nazarov, 
G.~I.~Olshanski and
V.~O.~Tarasov  
for reading the manuscript and valuable comments.

\section{The Yangian for the general linear Lie algebra}\label{sec:algstra}
\setcounter{equation}{0}

As we pointed out in the Introduction, the discovery of the Yangians was motivated
by the quantum inverse scattering theory.
It is possible, however, to ``observe"
the Yangian defining relations from purely algebraic viewpoint.
We start by showing that they are satisfied by certain natural
elements of the universal enveloping algebra $\U(\gl_n)$.
Then we show that these relations can be written
in a matrix form which provides a starting point for 
special algebraic techniques
to study the Yangian structure.

\subsection{Algebraic motivations and definitions}\label{subsec:mot}

Consider the general linear Lie algebra $\gl_n$ with its standard
basis $E_{ij}$, $i,j=1,\dots,n$. The commutation relations
are given by
\ben
[E_{ij},E_{kl}]=\delta_{kj}E_{il}-\delta_{il}E_{kj},
\end{equation*}
where $\delta_{ij}$ is the Kronecker delta; see e.g. \cite{h:il}.
Introduce the $n\times n$-matrix $E$ whose $ij$-th entry is $E_{ij}$.
The traces of powers of the matrix $E$ 
\ben
g_s=\tr E^{\tss s},\qquad s=1,2,\dots
\end{equation*}
are central elements
of the universal enveloping algebra
$\U(\gl_n)$ known as the {\it Gelfand invariants\/}; see \cite{g:tc}.
Moreover, the first $n$ of them are algebraically independent
and generate the center. A proof of the centrality of the $g_s$ is easily
deduced from the following relations in the enveloping algebra
\beql{commer}
[E_{ij}, (E^{\tss s})_{kl}]=\delta_{kj}(E^{\tss s})_{il}-\delta_{il}(E^{\tss s})_{kj}.
\end{equation}
One could wonder whether a more general closed formula
exists for the commutator of the matrix elements of 
the powers of $E$. The answer to this question turns out to be
affirmative and the following 
generalization of \eqref{commer} can be verified by induction:
\ben
[(E^{\tss r+1})_{ij}, (E^{\tss s})_{kl}]-[(E^{\tss r})_{ij}, (E^{\tss s+1})_{kl}]=
(E^{\tss r})_{kj}(E^{\tss s})_{il}-(E^{\tss s})_{kj}(E^{\tss r})_{il},
\end{equation*}
where $r,s\geq 0$ and $E^{\tss 0}=1$ is the identity matrix.
We can axiomatize these relations by introducing the following definition.

\bde\label{def:yangian}	The {\em Yangian for\/} $\gl_n$
is a unital associative algebra
with countably many generators $t_{ij}^{(1)},\ t_{ij}^{(2)},\dots$ where
$1\leq i,j\leq n$, and the defining relations
\beql{defyang}
[t^{(r+1)}_{ij}, t^{(s)}_{kl}]-[t^{(r)}_{ij}, t^{(s+1)}_{kl}]=
t^{(r)}_{kj} t^{(s)}_{il}-t^{(s)}_{kj} t^{(r)}_{il},
\end{equation}
where $r,s\geq 0$ and  $t^{(0)}_{ij}=\delta_{ij}$.
This algebra is denoted by
$\Y(\gl_n)$, or $\Y(n)$ for brevity. 
\ede
Introducing the generating
series,
\ben
t_{ij}(u) = \delta_{ij} + t^{(1)}_{ij} u^{-1} + t^{(2)}_{ij}u^{-2} +
\cdots\in\Y(n)[[u^{-1}]],
\end{equation*}
we can write \eqref{defyang} in the form
\beql{defrel}
(u-v)\ts [t_{ij}(u),t_{kl}(v)]=t_{kj}(u)t_{il}(v)-t_{kj}(v)t_{il}(u).
\end{equation}
Divide both sides by $u-v$ and use the formal expansion
\ben
(u-v)^{-1}=\sum_{r=0}^{\infty} u^{-r-1}v^r
\end{equation*}
to get an equivalent 
form of the defining relations
\beql{defequiv}
[t^{(r)}_{ij}, t^{(s)}_{kl}] =\sum_{a=1}^{\min(r,s)}
\big(t^{(a-1)}_{kj} t^{(r+s-a)}_{il}-t^{(r+s-a)}_{kj} t^{(a-1)}_{il}\big).
\end{equation}
The previous discussion
implies that the algebra $\Y(n)$ is nontrivial, as the mapping
\beql{mapyu}
t_{ij}^{(r)}\mapsto (E^{\tss r})_{ij}
\end{equation}
defines an algebra homomorphism $\Y(n)\to\U(\gl_n)$.

Alternatively, the generators of the Yangian can be realized
as the ``Capelli minors".
Keeping the notation $E$ for the matrix of the basis elements of $\gl_n$
introduce the {\it Capelli determinant\/}
\beql{capdetu}
\det (1+Eu^{-1})=\sum_{\sigma\in\Sym_n}\sgn\sigma\cdot (1+Eu^{-1})_{\sigma(1),1}\cdots
(1+E(u-n+1)^{-1})_{\sigma(n),n}.
\end{equation}
When multiplied by $u(u-1)\cdots (u-n+1)$ this determinant becomes a polynomial
in $u$ whose coefficients constitute another family of
algebraically independent generators of the center of $\U(\gl_n)$.	
The value of this polynomial at $u=n-1$ is a distinguished central element
whose image in a natural representation of $\gl_n$ by differential operators
is given by the celebrated {\it Capelli identity\/} \cite{c:zc}; see also \cite{hu:ci}.
For a positive integer $m\leq n$ introduce
the subsets of indices
${\mathcal B}_i=\{i,m+1,m+2,\dots,n\}$ and for any $1\leq i,j\leq m$ 
consider the Capelli minor
$\det (1+Eu^{-1})_{{\mathcal B}_i{\mathcal B}_j}$ defined as in 
\eqref{capdetu}, whose rows and columns are respectively enumerated
by 	${\mathcal B}_i$ and ${\mathcal B}_j$. These minors turn out to
satisfy the Yangian defining relations, i.e.,
there is an algebra
homomorphism
\ben
\Y(m)\to\U(\gl_{n}),\qquad t_{ij}(u)\mapsto 
\det (1+Eu^{-1})_{{\mathcal B}_i{\mathcal B}_j}.
\end{equation*}

These two interpretations of the Yangian defining relations
(which will reappear in Sections~\ref{subsec:pbw}, 
\ref{subsec:qst} and \ref{subsec:cc}) indicate
a close relationship between the representation theory of
the algebra $\Y(n)$ 
and the conventional representation theory of the general linear Lie algebra.
Many applications of the Yangian
are based on the following simple observation.

\bpr\label{prop:eval}
The mapping
\beql{eval}
\pi : t_{ij}(u)\mapsto \delta_{ij}+E_{ij}u^{-1}
\end{equation}
defines an algebra epimorphism $\Y(n)\to\U(\gl_n)$. Moreover,
\ben
E_{ij}\mapsto t_{ij}^{(1)}
\end{equation*}
is an embedding $\U(\gl_n)\hookrightarrow \Y(n)$.
\epr

In particular, any $\gl_n$-module can be extended to the algebra $\Y(n)$ via
\eqref{eval}. This plays an important role in the Yangian
representation theory.

\subsection{A matrix form of the defining relations}\label{subsec:matrix}

Introduce the $n\times n$ matrix $T(u)$ whose $ij$-th entry is the series 
$t_{ij}(u)$. It is convenient to regard it as an element
of the algebra $\Y(n)[[u^{-1}]]\ot \End \CC^n$ given by
\beql{tmatrix}
T(u)=\sum_{i,j=1}^n t_{ij}(u)\ot e_{ij},
\end{equation}
where the $e_{ij}$ denote the standard matrix 
units. For any positive integer $m$ we shall be using the algebras of the form
\beql{multitp}
\Y(n)[[u^{-1}]]\ot \End \CC^n\ot \cdots\ot \End \CC^n,
\end{equation}
with $m$ copies of $\End \CC^n$. For any $a\in\{1,\dots,m\}$ we denote by $T_a(u)$
the matrix $T(u)$ which acts on the $a$-th copy of $\End \CC^n$. That is,
$T_a(u)$ is an element of the algebra \eqref{multitp} of the form
\beql{tku}
T_a(u)=\sum_{i,j=1}^n t_{ij}(u)\ot 1\ot \cdots \ot 1\ot e_{ij}\ot 1\ot \cdots\ot 1,
\end{equation}
where the $e_{ij}$ belong to the $a$-th copy of $\End \CC^n$ and 
$1$ is the identity matrix. Similarly, if
\ben
C=\sum_{i,j,k,l=1}^n c^{}_{ijkl}\ts e_{ij}\ot e_{kl}\in \End \CC^n\ot\End \CC^n,
\end{equation*}
then for distinct indices $a,b\in\{1,\dots,m\}$ we introduce
the element $C_{ab}$ of the algebra \eqref{multitp} by
\beql{ars}
C_{ab}=\sum_{i,j,k,l=1}^n c^{}_{ijkl}\ts 
1\ot 1\ot \cdots \ot 1\ot e_{ij}\ot 1\ot \cdots\ot 1\ot e_{kl}\ot 1\ot \cdots\ot 1,
\end{equation}
where the  $e_{ij}$ and $e_{kl}$ belong to the $a$-th and $b$-th copies of 
$\End \CC^n$, respectively.
Consider now the permutation operator
\ben
P=\sum_{i,j=1}^n e_{ij}\ot e_{ji}\in \End \CC^n\ot\End \CC^n.
\end{equation*}
The rational function
\beql{rmatrix}
R(u)=1-Pu^{-1}
\end{equation}
with values in $\End \CC^n\ot\End \CC^n$ is called 
the {\it Yang\/} $R$-{\it matrix\/}. (Here and below we write $1$ instead of 
$1\ot 1$ for brevity).
An easy calculation in the group algebra $\CC[\Sym_3]$ shows that
the following identity holds in the algebra $(\End \CC^n)^{\ot 3}$
\beql{ybe}
R_{12}(u)R_{13}(u+v)R_{23}(v)=R_{23}(v)R_{13}(u+v)R_{12}(u).
\end{equation}
This is known as the {\it Yang--Baxter equation\/}. The Yang $R$-matrix
is its simplest nontrivial solution. In the following we regard
$T_1(u)$ and $T_2(u)$ as elements of the algebra \eqref{multitp} with $m=2$.

\bpr\label{prop:ternary}
The defining relations \eqref{defyang} can be written in the equivalent 
form
\beql{ternary}
R(u-v)T_1(u)T_2(v)=T_2(v)T_1(u)R(u-v).
\end{equation}
\epr

The relation \eqref{ternary} is known 
as the {\it ternary\/} or $RTT$ {\it relation\/}.

\subsection{Automorphisms and anti-automorphisms}\label{subsec:auto}

Consider an arbitrary formal series	which begins with 1,
\ben
f(u)=1+f_1u^{-1}+f_2 u^{-2}+\cdots\in \CC[[u^{-1}]].
\end{equation*}
Also, let $a\in\CC$ be an arbitrary element and $B$ an arbitrary nondegenerate
complex matrix.

\bpr\label{prop:auto}
Each of the mappings
\begin{align}
\label{mult}
T(u)&\mapsto f(u) T(u),\\
\label{shift}
T(u)&\mapsto T(u+a), \\
T(u)&\mapsto BT(u)B^{-1}
\non
\end{align}
defines an automorphism of $\Y(n)$.
\epr

The matrix $T(u)$ can be regarded as a formal series in $u^{-1}$
whose coefficients are matrices over $\Y(n)$. Since this series
begins with the identity matrix, $T(u)$ is invertible
and we denote by $T^{-1}(u)$ the inverse element.
Also, denote by $T^{\tss t}(u)$ the transposed matrix for $T(u)$.

\bpr\label{prop:anti}
Each of the mappings
\begin{align}
T(u)&\mapsto T(-u),
\non\\
T(u)&\mapsto T^{\tss t}(u), 
\non\\
T(u)&\mapsto T^{-1}(u)
\non
\end{align}
defines an anti-automorphism of $\Y(n)$.
\epr

\subsection{Poincar\'e--Birkhoff--Witt theorem}\label{subsec:pbw}

\bth\label{thm:pbw}
Given an arbitrary linear order on the set of
the generators
$t^{(r)}_{ij}$, any element of the algebra $\Y(n)$ is uniquely written
as a linear combination of ordered monomials in the generators.
\eth

\Outline There are two natural ascending filtrations on 
the algebra $\Y(n)$. Here we use the one
defined by
\ben
\deg_1 t_{ij}^{(r)}=r.
\end{equation*}
(The other filtration will be used in Section~\ref{subsec:center}). It is immediate from
the defining relations \eqref{defequiv} that the corresponding graded
algebra $\gr_1 \Y(n)$ is commutative. Denote by $\overline{t}_{ij}^{\ts(r)}$
the image of $t_{ij}^{(r)}$ in the $r$-th component of $\gr_1 \Y(n)$.
It will be sufficient to show that the elements $\overline{t}_{ij}^{\ts(r)}$
are algebraically independent.

The composition of the automorphism $T(u)\mapsto T^{-1}(-u)$ of $\Y(n)$
and the homomorphism \eqref{eval} yields another homomorphism
$\Y(n)\to\U(\gl_n)$ such that
\beql{homtinv}
T(u)\mapsto (1-Eu^{-1})^{-1}.
\end{equation}
The image of the generator $t_{ij}^{(r)}$ is given by
\eqref{mapyu}.
For any nonnegative integer $m$ consider the Lie algebra $\gl_{n+m}$
and now let $E$ denote the corresponding matrix formed by its basis elements
$E_{ij}$. Then formula \eqref{mapyu} still defines a homomorphism
\beql{pim}
\Y(n)\to\U(\gl_{n+m}). 
\end{equation}
Consider the canonical filtration on the
universal enveloping algebra $\U(\gl_{n+m})$ and observe that
the homomorphism \eqref{pim} is filtration-preserving. So, it defines
a homomorphism of the corresponding graded algebras
\beql{pimbar}
\gr_1 \Y(n)\to \Sr(\gl_{n+m}),
\end{equation}
where  $\Sr(\gl_{n+m})$ is the symmetric algebra of $\gl_{n+m}$.
One can show that for any finite family of elements $\overline{t}_{ij}^{\ts(r)}$
there exists a sufficiently large parameter $m$ such that their images
under \eqref{pimbar} 
are algebraically independent. 
\endproof

\bco\label{cor:polalg}
$\gr_1 \Y(n)$ is the algebra of polynomials in the variables 
$\overline{t}_{ij}^{\ts(r)}$. 
\eco

\subsection{Hopf algebra structure}\label{subsec:hopf}

A {\it Hopf algebra\/} $A$ (over $\CC$) is an associative algebra equipped
with a {\it coproduct\/} (or {\it comultiplication}) $\Delta:A\to A\ot A$, 
an {\it antipode\/} $\Sr:A\to A$ and a {\it counit\/}
$\epsilon:A\to\CC$ such that $\Delta$ and $\epsilon$ are algebra homomorphisms,
$\Sr$ is an anti-automorphism and some other axioms are satisfied.
These can be found in any textbook on the subject.
In \cite[Chapter~4]{cp:gq} the Hopf algebra axioms are discussed in the context of
quantum groups.

\bth\label{thm:hopf}
The Yangian $\Y(n)$ is a Hopf algebra with the coproduct
\ben
\Delta: t_{ij}(u)\mapsto \sum_{a=1}^n t_{ia}(u)\ot t_{aj}(u),
\end{equation*}
the antipode
\ben
\Sr: T(u)\mapsto T^{-1}(u),
\end{equation*}
and the counit
\ben
\epsilon: T(u)\mapsto 1.
\end{equation*}
\eth

\Proof We only verify the most nontrivial axiom that $\Delta:\Y(n)\to\Y(n)\ot\Y(n)$
is an algebra homomorphism. The remaining axioms follow
directly from the definitions. We slightly generalize the notation used in 
Section~\ref{subsec:matrix}. Let $p$ and $m$ be positive integers.
Introduce the algebra
\ben
\big(\Y(n)[[u^{-1}]]\big)^{\otimes p}\otimes \big(\End\CC^n\big)^{\otimes m}
\end{equation*}
and for all $a\in\{1,\dots,m\}$ and $b\in\{1,\dots,p\}$ consider
its elements
\ben
T_{[b]a}(u)=\sum_{i,j=1}^n (1^{\otimes b-1}\otimes t_{ij}(u)
\otimes 1^{\otimes p-b})\otimes(1^{\otimes a-1}\otimes e_{ij}
\otimes 1^{\otimes m-a}).
\end{equation*}
The definition of $\Delta$ can now be written in a matrix form,
\ben
\Delta: T(u)\mapsto T_{[1]}(u)\ts T_{[2]}(u),
\end{equation*}
where $T_{[b]}(u)$ is an abbreviation for $T_{[b]1}(u)$.
We need to show that
that $\Delta (T(u))$ satisfies the ternary relation \eqref{ternary}, i.e.,
\ben
R(u-v)T_{[1]1}(u)T_{[2]1}(u)T_{[1]2}(v)T_{[2]2}(v)=
T_{[1]2}(v)T_{[2]2}(v)T_{[1]1}(u)T_{[2]1}(u)R(u-v).
\end{equation*}
However, this is implied by the ternary relation \eqref{ternary}
and the observation that $T_{[2]1}(u)$ and
$T_{[1]2}(v)$, as well as $T_{[1]1}(u)$ and $T_{[2]2}(v)$, commute.
\endproof

\subsection{Quantum determinant and quantum minors}\label{subsec:qdet}

Consider the rational function $R(u_1,\dots,u_m)$ with values in the algebra
$(\End\CC^n)^{\ot m}$ defined by
\beql{prodrma}
R(u_1,\dots,u_m)=(R_{m-1,m})(R_{m-2,m}R_{m-2,m-1}) \cdots (R_{1m}
\cdots
R_{12}),
\end{equation}
where we abbreviate $R_{ij}=R_{ij}(u_i-u_j)$. Applying the Yang--Baxter equation
\eqref{ybe} and the fact that $R_{ij}$ and $R_{kl}$ commute if the indices are distinct,
we can write \eqref{prodrma} in a different form. In particular,
\ben
R(u_1,\dots,u_m)=(R_{12}\cdots R_{1m})\cdots (R_{m-2,m-1}R_{m-2,m})
(R_{m-1,m}).
\end{equation*}
As before, we use the notation $T_a(u_a)$
for the matrix $T(u_a)$ of the Yangian generators
corresponding to the $a$-th copy of $\End\CC^n$.

\bpr\label{prop:fundam}
We have the relation
\ben
R(u_1,\dots, u_m)\, T_1 (u_1) \cdots T_m (u_m) = T_m (u_m) \cdots T_1
(u_1)\, R(u_1,\dots, u_m).
\end{equation*}
\epr

We let the $e_i$, $i=1,\dots,n$ denote the canonical basis of $\CC^n$, and
$A_m$ the antisymmetrizer in $(\CC^n)^{\ot m}$ given by
\beql{antisym}
A_m(e_{i_1}\ot\cdots\ot e_{i_m})=\sum_{\sigma\in\Sym_m}\sgn\sigma\cdot
e_{i_{\sigma(1)}}\ot\cdots	\ot e_{i_{\sigma(m)}}.
\end{equation}
Note that this operator satisfies $A_m^2=m!\ts A_m$.

\bpr\label{prop:rmanti}
If $u_i-u_{i+1}=1$ for all $i=1,\dots,m-1$ then
\ben
R(u_1,\dots, u_m)=A_m.
\end{equation*}
\epr

By Propositions~\ref{prop:fundam}
and \ref{prop:rmanti} we have the identity
\beql{att}
A_m\, T_1 (u) \cdots T_m (u-m+1) = T_m (u-m+1)\cdots T_1 (u) \, A_m.
\end{equation}
Suppose now that $m=n$. Then the antisymmetrizer is a one-dimensional operator
in $(\CC^n)^{\ot n}$. 
Therefore, the element \eqref{att} equals $A_n$ times a scalar series with coefficients
in $\Y(n)$ which prompts the following definition.

\bde\label{def:qdet}
The {\em quantum determinant\/} is the formal series
\beql{qdet}
\qdet T(u)=1+d_1 u^{-1}+d_2 u^{-2}+\cdots\in\Y(n)[[u^{-1}]]
\end{equation}
such that the element \eqref{att} {\em(}with $m=n${\em)} equals $A_n\ts\qdet T(u)$.
\ede

\bpr\label{prop:expqdet}
For any permutation $\rho\in\Sym_n$ we have
\begin{align}\label{qdet1}
\qdet T(u)&=\sgn \rho \sum_{\sigma\in\Sym_n}\sgn\sigma\cdot
t_{\sigma(1),\rho(1)}(u)\cdots t_{\sigma(n),\rho(n)}(u-n+1)\\
\label{qdet2}
{}&=\sgn \rho \sum_{\sigma\in\Sym_n}\sgn\sigma\cdot
t_{\rho(1),\sigma(1)}(u-n+1)\cdots t_{\rho(n),\sigma(n)}(u).
\end{align}
\epr

\bex\label{ex:qdet}
In the case $n=2$ we have
\ben
\bal
\qdet T(u)&=t_{11}(u)\ts t_{22}(u-1)-t_{21}(u)\ts t_{12}(u-1)\\
{} &=t_{22}(u)\ts t_{11}(u-1)-t_{12}(u)\ts t_{21}(u-1)\\
{} &=t_{11}(u-1)\ts t_{22}(u)-t_{12}(u-1)\ts t_{21}(u)\\
{} &=t_{22}(u-1)\ts t_{11}(u)-t_{21}(u-1)\ts t_{12}(u).
\eal
\end{equation*}
\eex

More generally, assuming that $m\leq n$ is arbitrary, we can define
$m\times m$ quantum minors as the matrix elements of the operator \eqref{att}.
Namely, the operator \eqref{att} can be written as
\ben
\sum {t\ts}^{c_1\cdots\ts c_m}_{d_1\cdots\ts d_m}(u)\ot e_{c_1d_1}\ot \cdots
\ot e_{c_md_m},
\end{equation*}
summed over the indices $c_i,d_i\in\{1,\dots,n\}$, where
${t\ts}^{c_1\cdots\ts c_m}_{d_1\cdots\ts d_m}(u)\in \Y(n)[[u^{-1}]]$.
We call these elements the {\it quantum	minors\/} of the matrix $T(u)$.
The following formulas are obvious generalizations of 
\eqref{qdet1} and \eqref{qdet2},
\ben
\bal
{t\ts}^{c_1\cdots\ts c_m}_{d_1\cdots\ts d_m}(u)&=
\sum_{\sigma\in \Sym_m} \sgn\sigma\cdot t_{c_{\sigma(1)}d_1}(u)\cdots
t_{c_{\sigma(m)}d_m}(u-m+1)\\
{}&=
\sum_{\sigma\in \Sym_m} \sgn\sigma\cdot t_{c_1d_{\sigma(1)}}(u-m+1)\cdots
t_{c_md_{\sigma(m)}}(u).
\eal
\end{equation*}
It is clear from the definition that the quantum minors are skew-symmetric
with respect to permutations of the upper indices and of the lower indices.

\bpr\label{prop:coprqm}
The images of the quantum minors under the coproduct are given by
\beql{coprqm}
\Delta ({t\ts}^{c_1\cdots\ts c_m}_{d_1\cdots\ts d_m}(u))
=\sum_{a_1<\cdots < a_m}{t\ts}^{c_1\cdots\ts c_m}_{a_1\cdots\ts a_m}(u)
\ot {t\ts}^{a_1\cdots\ts a_m}_{d_1\cdots\ts d_m}(u),
\end{equation}
summed over all subsets of indices $\{a_1,\dots,a_m\}$ from $\{1,\dots,n\}$.
\epr

\Proof Using the notation of Section~\ref{subsec:hopf} we can write the image of
the left hand side of \eqref{att} under the coproduct $\Delta$ as
\ben
A_m\, T_{[1]1} (u)T_{[2]1} (u) \cdots T_{[1]m} (u-m+1)T_{[2]m} (u-m+1).
\end{equation*}
Since $m!\ts A_m=A_m^2$, by \eqref{att} this coincides with
\ben
\frac{1}{m!}\ts
A_m\ts T_{[1]1} (u)\cdots T_{[1]m} (u-m+1)
\ts A_m\ts T_{[2]1} (u) \cdots T_{[2]m} (u-m+1).
\end{equation*}
Taking here the matrix elements and using the skew-symmetry of the quantum minors
we come to \eqref{coprqm}. \endproof

\bco\label{cor:delqdet}
We have
\ben
\Delta: \qdet T(u)\mapsto \qdet T(u)\ot\qdet T(u).
\end{equation*}
\eco

\subsection{The center of $\Y(n)$}\label{subsec:center}

\bpr\label{prop:qminorel}
We have the relations
\ben
[{t}^{}_{ab}(u),
{t\ts}^{c_1\cdots\ts c_m}_{d_1\cdots\ts d_m}(v)]=
\frac{1}
{u-v}
\left(\sum_{i=1}^m
{t}^{}_{c_ib}(u)
{t\ts}^{c_1\cdots\ts a\ts\cdots\ts c_m}_{d_1\ \cdots\ \  d_m}(v)
-\sum_{i=1}^m{t\ts}^{c_1\ \cdots\ \  c_m}_{d_1\cdots\ts 
b\ts\cdots\ts d_m}(v)
{t}^{}_{ad_i}(u)\right),
\end{equation*}
where the indices $a$ and $b$ in the quantum minors
replace $c_i$ and $d_i$, respectively.
\epr

\bco\label{cor:qmincent} For any indices $i,j$ we have
\ben
[{t}^{}_{c_id_j}(u),
{t\ts}^{c_1\cdots\ts c_m}_{d_1\cdots\ts d_m}(v)]=0.	
\end{equation*}
\eco

Recall the elements $d_i\in\Y(n)$ are defined by \eqref{qdet}.

\bth\label{thm:qdcenter}
The coefficients $d_1,d_2,\dots$ of the series $\qdet T(u)$
belong to the center of the algebra $\Y(n)$. Moreover, these elements are algebraically
independent and generate the center of $\Y(n)$.
\eth

\Outline The first claim follows from Corollary~\ref{cor:qmincent}.
To prove the second claim introduce a filtration on $\Y(n)$ by setting
\ben
\deg_2 t_{ij}^{(r)}=r-1.
\end{equation*}
The the corresponding graded algebra $\gr_2 \Y(n)$ is isomorphic to the
universal enveloping algebra $\U(\gl_n[x])$ where $\gl_n[x]$ is the
Lie algebra of polynomials in an indeterminate $x$ with coefficients in $\gl_n$.
Indeed, denote by $\overline{t}_{ij}^{\ts(r)}$ the image of $t_{ij}^{(r)}$
in the $(r-1)$-th component of $\gr_2 \Y(n)$. Then 
$E_{ij}\ts x^{r-1}\mapsto \overline{t}_{ij}^{\ts(r)}$ is
an algebra homomorphism $\U(\gl_n[x])\to \gr_2 \Y(n)$. Its kernel
is trivial by Theorem~\ref{thm:pbw}.

We observe now from \eqref{qdet1} (with $\rho=1$) that the coefficient $d_r$
of $\qdet T(u)$ has the form
\ben
d_r=t_{11}^{(r)}+\cdots+ t_{nn}^{(r)} +\ \ \text{terms of degree}\ \ <r-1.
\end{equation*}
This implies that the elements $d_r$, $r\geq 1$ are algebraically independent.
Furthermore, the image of $d_r$ in the $(r-1)$-th component of $\gr_2 \Y(n)$
coincides with $Z\ts x^{r-1}$ where $Z=E_{11}+\cdots+E_{nn}$. It remains to note
that the elements $Z\ts x^{r-1}$ with $r\geq 1$ generate the center of
$\U(\gl_n[x])$. The latter follows from the fact that the center of $\U(\sll_n[x])$
is trivial \cite{ci:hc}; see also \cite[Proposition~2.12]{mno:yc}. \endproof

\subsection{The Yangian for the special linear Lie algebra}\label{subsec:sln}

For any simple Lie algebra $\mathfrak{a}$ over $\CC$ the 
corresponding Yangian $\Y(\mathfrak{a})$
is a deformation of the universal enveloping algebra $\U(\mathfrak{a}[x])$
in the class of Hopf algebras. Two different
presentations of $\Y(\mathfrak{a})$ are given in \cite{d:ha} and \cite{d:nr}.
The type $A$ Yangian $\Y(\sll_n)$ can also be realized as a Hopf subalgebra
of $\Y(n)$, as well as a quotient of $\Y(n)$.

For any series $f(u)\in 1+u^{-1}\CC[[u^{-1}]]$ denote by $\mu_f$
the automorphism \eqref{mult} of $\Y(n)$.

\bde\label{def:sln}
The {\em Yangian for\/} $\sll_n$ is the subalgebra $\Y(\sll_n)$ of $\ts\Y(n)$
which consists of the elements stable under all automorphisms
$\mu_f$.
\ede

We let $\Z(n)$ denote the center of $\Y(n)$.

\bth\label{thm:decoyn}
The subalgebra $\Y(\sll_n)$ is a Hopf algebra
whose coproduct, antipode and counit are obtained by restricting those
from $\Y(n)$. Moreover,	 $\Y(n)$ is isomorphic to the tensor product
of its subalgebras
\beql{decoyn}
\Y(n)=\Z(n)\ot 	\Y(\sll_n).
\end{equation}
\eth

\Outline It is easy to verify
that there exists a unique formal power series
$\wt d(u)$ in $u^{-1}$ with coefficients in $\Z(n)$ 
which begins with 1 and satisfies
\ben
\wt d(u)\ts\wt d(u-1)\cdots \wt d(u-n+1)=\qdet T(u).
\end{equation*}
Then by Proposition~\ref{prop:expqdet}, the image of $\wt d(u)$ under $\mu_f$
is given by
\beql{dtimuf}
\mu_f: 	\wt d(u)\mapsto f(u)\ts \wt d(u).
\end{equation}
This implies that all coefficients of the series
\ben
\tau_{ij}(u)=\wt d(u)^{-1}\ts t_{ij}(u)
\end{equation*}
belong to the subalgebra $\Y(\sll_n)$. In fact, they generate
this subalgebra. Furthermore, 
the coefficients of the series $\wt d(u)$
are algebraically independent over $\Y(\sll_n)$ which gives
\eqref{decoyn}.
Corollary~\ref{cor:delqdet} implies that 
\ben
\Delta: \wt d(u)\mapsto \wt d(u)\ot\wt d(u)
\end{equation*}
and so the image of  $\Y(\sll_n)$
under the coproduct is contained in $\Y(\sll_n)\ot\Y(\sll_n)$.
Using Definition~\ref{def:qdet}, we find that
the image of $\qdet T(u)$ under the antipode $\Sr$ is $(\qdet T(u))^{-1}$
and so,
\ben
\Sr : \wt d(u)^{-1}\ts T(u)\mapsto \wt d(u)\ts T^{-1}(u).
\end{equation*}
Due to \eqref{dtimuf}, $\Y(\sll_n)$ is stable under $\Sr$.
\endproof

\bco\label{cor:quoti}
The algebra $\Y(\sll_n)$ is isomorphic to the quotient of \ts $\Y(n)$ 
by the ideal generated by the center, i.e.,
\ben
\Y(\sll_n)\cong \Y(n)/(\qdet T(u)=1).
\end{equation*}
\eco

\subsection{Two more realizations of the Hopf algebra $\Y(\sll_2)$}
\label{subsec:3real}

\bde\label{def:adef1}
Let $\mathcal A$ denote the associative unital algebra
with six generators $e,f,h,J(e),J(f),J(h)$ and the defining relations
\ben
\begin{aligned}[]
[e,f]&=h,\qquad [h,e]=2e,\qquad [h,f]=-2f,\\
[x,J(y)]&=J([x,y]),\qquad J(ax)=a\tss J(x),
\end{aligned}
\end{equation*}
where $x,y\in\{e,f,h\}$, $a\in\CC$, and 
\ben
\big[[J(e),J(f)],J(h)\big]=\big(J(e)f-e\ts J(f)\big)\tss h.	 
\end{equation*}
The Hopf algebra structure on $\mathcal A$ is defined by
\ben
\bal
\Delta:x&\mapsto x\ot 1+1\ot x,\\
\Delta:J(x)&\mapsto J(x)\ot 1+1\ot J(x)+\frac12[x\ot 1,t],\\
\Sr:x&\mapsto -x,\qquad J(x)\mapsto -J(x)+x,\\
\epsilon:x&\mapsto 0,\qquad J(x)\mapsto 0,
\eal
\end{equation*}
where $t=e\ot f+f\ot e+\frac12 h\ot h$.
\ede

This definition generalizes to any simple Lie algebra $\agot$. 
Given a basis $e_1,\dots,e_n$ of $\agot$,
the corresponding Yangian $\Y(\agot)$ is a Hopf algebra
generated by $2n$ elements $e_i$, $J(e_i)$, $i=1,\dots,n$ 
as originally defined by Drinfeld~\cite{d:ha}.

\bth\label{thm:isoma1}
The mapping 
\ben
e\mapsto t_{12}^{(1)},\qquad f\mapsto t_{21}^{(1)},\qquad h\mapsto 
t_{11}^{(1)}-t_{22}^{(1)},
\end{equation*}
\ben
\bal
J(e)&\mapsto t_{12}^{(2)}-\frac12(t_{11}^{(1)}+t_{22}^{(1)}-1)\ts t_{12}^{(1)},\\
J(f)&\mapsto t_{21}^{(2)}-\frac12(t_{11}^{(1)}+t_{22}^{(1)}-1)\ts t_{21}^{(1)},\\
J(h)&\mapsto t_{11}^{(2)}-t_{22}^{(2)}-
\frac12(t_{11}^{(1)}+t_{22}^{(1)}-1)(t_{11}^{(1)}-t_{22}^{(1)}).
\eal
\end{equation*}
defines a Hopf algebra isomorphism $\mathcal A\to\Y(\sll_2)$.
\eth

\bde\label{def:ydefgf2}
Let $\mathcal B$ be the associative algebra
with generators $e_k,f_k,h_k$ where $k=0,1,2,\dots$ and the defining relations
given in terms of the generating series
\ben
\bal
e(u)&=\sum_{k=0}^{\infty}e_k\ts u^{-k-1},\\
f(u)&=\sum_{k=0}^{\infty}f_k\ts u^{-k-1},\\
h(u)&=1+\sum_{k=0}^{\infty}h_k\ts u^{-k-1}
\eal
\end{equation*}
as follows
\ben
[h(u),h(v)]=0,\qquad [e(u),f(v)]=-\frac{h(u)-h(v)}{u-v},
\end{equation*}
\ben
\bal[]
[e(u),e(v)]&=-\frac{\big(e(u)-e(v)\big)^2}{u-v},\\
[f(u),f(v)]&=\frac{\big(f(u)-f(v)\big)^2}{u-v},\\
[h(u),e(v)]&=-\frac{\big\{h(u),e(u)-e(v)\big\}}{u-v},\\
[h(u),f(v)]&=\frac{\big\{h(u),f(u)-f(v)\big\}}{u-v},
\eal
\end{equation*}
where we have used the notation $\{a,b\}=ab+ba$.
The Hopf algebra structure on $\mathcal B$ is defined by
the coproduct
\ben
\bal
\Delta: e(u)&\mapsto e(u)\ot 1+
\sum_{k=0}^{\infty}(-1)^k f(u+1)^k\ts h(u)\ot
e(u)^{k+1},\\
\Delta: f(u)&\mapsto 1\ot f(u) +
\sum_{k=0}^{\infty}(-1)^k f(u)^{k+1}\ot h(u)\ts 
e(u+1)^{k},\\
\Delta: h(u)&\mapsto
\sum_{k=0}^{\infty}(-1)^k(k+1)f(u+1)^{k}\ts h(u)\ot 
h(u)\ts e(u+1)^{k},
\eal
\end{equation*}
the antipode
\ben
\bal
\Sr: e(u)&\mapsto -\big(h(u)+f(u+1)e(u)\big)^{-1}\ts e(u),\\
\Sr: f(u)&\mapsto - f(u)\ts\big(h(u)+f(u)e(u+1)\big)^{-1},\\
\Sr: h(u)&\mapsto \big(h(u)+f(u+1)e(u)\big)^{-1}
\Big(1-f(u+1)\big(h(u)+f(u+1)e(u)\big)^{-1}e(u)\Big),
\eal
\end{equation*}
and	the counit
\ben
\epsilon: e(u)\mapsto 0,\qquad	f(u)\mapsto 0,\qquad h(u)\mapsto 1.
\end{equation*}
\ede

\noindent
Explicitly, the defining relations of ${\mathcal B}$
can be written in the form
\ben
[h_k,h_l]=0,\qquad [e_k,f_l]=h_{k+l},
\qquad [h_0,e_k]=2e_k,\qquad [h_0,f_k]=-2f_k,
\end{equation*}
\ben
\bal[]
[e_{k+1},e_l]-[e_{k},e_{l+1}]&=e_{k}e_l+e_{l}e_k,\\
[f_{k+1},f_l]-[f_{k},f_{l+1}]&=-f_{k}f_l-f_{l}f_k,\\
[h_{k+1},e_l]-[h_{k},e_{l+1}]&=h_{k}e_l+e_{l}h_k,\\
[h_{k+1},f_l]-[h_{k},f_{l+1}]&=-h_{k}f_l-f_{l}h_k.
\eal
\end{equation*}

Such a realization exists for an arbitrary Yangian $\Y(\agot)$. Some authors
call it the {\it new realization\/} following the title of Drinfeld's paper~\cite{d:nr}
where it was introduced. This presentation of the Yangian is most convenient
to describe its finite-dimensional irreducible representations \cite{d:nr}.
However, in the case of an arbitrary
simple Lie algebra $\agot$ no explicit formulas for 
the coproduct and antipode are known.

\bth\label{thm:isomb2}
The mapping 
\ben
\bal
e(u)&\mapsto t_{22}(u)^{-1}t_{12}(u),\\
f(u)&\mapsto t_{21}(u)\ts t_{22}(u)^{-1},\\
h(u)&\mapsto t_{11}(u)\ts t_{22}(u)^{-1}-t_{21}(u)
\ts t_{22}(u)^{-1}t_{12}(u)\ts t_{22}(u)^{-1},
\eal
\end{equation*}
defines a Hopf algebra isomorphism ${\mathcal B}\to\Y(\sll_2)$.
\eth

Combining the two above theorems we obtain a Hopf algebra isomorphism
${\mathcal A}\to{\mathcal B}$ given by
\ben
e\mapsto e_0,\qquad f\mapsto f_0,\qquad h\mapsto 
h_0,
\end{equation*}
and
\ben
\bal
J(e)&\mapsto e_1-\frac14(e_0\ts h_0+h_0\ts e_0),\\
J(f)&\mapsto f_1-\frac14(f_0\ts h_0+h_0\ts f_0),\\
J(h)&\mapsto h_1+\frac12(e_0\ts f_0+f_0\ts e_0-h_0^2).
\eal
\end{equation*}

As we have seen in the proof of Theorem~\ref{thm:qdcenter}, the graded
algebra $\gr_2 \Y(\sll_2)$ is isomorphic to the universal enveloping algebra
$\U(\sll_2[x])$. The images of the generators of the algebra ${\mathcal A}$
in the graded algebra clearly correspond to $e,f,h,ex,fx,hx$ while
the images of the generators $e_k,f_k,h_k$ of ${\mathcal B}$
correspond to $ex^k,fx^k,hx^k$.

\subsection{Quantum Liouville formula}\label{subsec:qlf}

Here we give another family of generators of the center of $\Y(n)$.
Introduce the series $z(u)$ with coefficients from $\Y(n)$ by the formula
\beql{qcont}
z(u)^{-1}=\frac{1}{n}\ts\tr\big(T(u)\ts T^{-1}(u-n)\big),
\end{equation}
so that
\ben
z(u)=1+z_{2}u^{-2}+z_3u^{-3}+\cdots, \qquad z_i\in\Y(n).
\end{equation*}

\bde\label{def:comatrix}
The {\em quantum comatrix\/} $\wh T(u)$ is defined by
\beql{comatrix}
\wh T(u)\ts T(u-n+1)=\qdet T(u).
\end{equation}
\ede

\bpr\label{prop:mateco}
The matrix elements $\wh t_{ij}(u)$ of the matrix $\wh T(u)$ are given by
\beql{mateco}
\wh t_{ij}(u)=(-1)^{i+j}
{t\ts}^{1\ts\cdots\ts
\wh j\ts\cdots\ts n}_{1\ts\cdots\ts\wh i\ts\cdots\ts n}(u),
\end{equation}
where the hats on the right hand side
indicate the indices to be omitted.	Moreover, we have the relation
\beql{cotrans}
\wh T^{\ts t}(u-1)\ts T^{\tss t}(u)=\qdet T(u).
\end{equation}
\epr

\Proof Using Definition~\ref{def:qdet} we derive from \eqref{comatrix}
\ben
A_n\ts T_1(u)\cdots T_{n-1}(u-n+2)=A_n\wh T_n(u).
\end{equation*}
Taking the matrix elements we come to \eqref{mateco}.  
Further, consider the automorphism $\varphi:T(u)\mapsto T^{\tss t}(-u)$ of $\Y(n)$; 
see Proposition~\ref{prop:anti}.
Using Proposition~\ref{prop:expqdet} we find
\ben
\varphi: \qdet T(u)\mapsto \qdet T(-u+n-1),\qquad 
\varphi: \wh T(u) \mapsto \wh T^{\ts t}(-u+n-2).
\end{equation*}
Now applying $\varphi$ to \eqref{comatrix} and replacing $-u+n-1$ by $u$
we get \eqref{cotrans}. \endproof

The following is a `quantum' analog of the classical Liouville formula;
see \cite[Remark~5.8]{mno:yc} for more comment.

\bth\label{thm:liouville}
We have the relation
\beql{liouville}
z(u)=\frac{\qdet T(u-1)}{\qdet T(u)}.
\end{equation}
\eth

\Proof From \eqref{qcont} and \eqref{comatrix} we find 
\ben
z(u)^{-1}=\frac{1}{n}\ts\tr\big(T(u)\ts \wh T(u-1)\ts (\qdet T(u-1))^{-1}\big).
\end{equation*}
Using the centrality of $\qdet T(u)$ and \eqref{cotrans} we get \eqref{liouville}.
\endproof

\bco\label{cor:zugencent}
The coefficients	$z_2,z_3,\dots$ of $z(u)$ are algebraically independent
generators of the center of $\Y(n)$.
\eco

\bpr\label{prop:sqantipode}
The square of the antipode $\Sr$ is the automorphism of $\Y(n)$ given by
\ben
\Sr^2: T(u)\mapsto z(u+n)\ts T(u+n).
\end{equation*}
In particular, $\qdet T(u)$ is stable under $\ts\Sr^2$.
\epr

\Outline  The series $z(u)$ can be defined equivalently in a way similar to
the quantum determinant; cf. Definition~\ref{def:qdet}.
Multiply the ternary relation \eqref{ternary}
by $T_2^{-1}(v)$ from both sides and take the transposition with respect to the second
copy of $\End\CC^n$. We obtain the relation
\beql{trater}
R^t(u-v) \ts\wt T_2(v)\ts T_1(u)  =	T_1(u)\ts \wt T_2(v)\ts R^t(u-v),
\end{equation}
where $\wt T(v)=\big(T^{-1}(v)\big)^t$ and
\ben
R^t(u)=1-{Q}{u^{-1}},\qquad Q=\sum_{i,j=1}^n e_{ij}\ot e_{ij}.
\end{equation*}
Observe that $Q$ is a one-dimensional operator satisfying $Q^2=n\ts Q$.
Therefore, $R^t(u)^{-1}$ has a simple pole at $u=n$ with residue $Q$.
Relation \eqref{trater} now implies
\beql{qresi}
Q\ts T_1(u)\ts \wt T_2(u-n)=\wt T_2(u-n) \ts T_1(u)\ts Q,
\end{equation}
and this element equals $Q\ts z(u)^{-1}$. It now suffices to apply 
the antipode $\Sr$
to both sides of the identity $T(u)\ts T^{-1}(u)=1$.
The second claim follows from Theorem~\ref{thm:liouville}.
\endproof

\subsection{Factorization of the quantum determinant}
\label{subsec:qdf}

\bde\label{def:quasidet}
Let $X$ be a square matrix over a ring with 1. Suppose that there exists
the inverse matrix $X^{-1}$ and its $ji$-th entry $(X^{-1})_{ji}$
is an invertible element of the ring. Then the $ij$-{\em th quasi-determinant
of} $X$ is defined by the formula 
\ben
|X|_{ij}=	\big((X^{-1})_{ji}\big)^{-1}.
\end{equation*}
\ede

For any $1\leq m\leq n$ denote by $T^{(m)}(u)$ the submatrix of $T(u)$
corresponding to the first $m$ rows and columns.

\bth\label{thm:qdfactor} The quantum determinant $\qdet T(u)$
admits the factorization in the algebra $\Y(n)[[u^{-1}]]$
\beql{qdfactor}
\qdet T(u)=t_{11}(u) \ts |T^{(2)}(u-1)|_{22}\cdots |T^{(n)}(u-n+1)|_{nn}.
\end{equation}
Moreover, the factors on the right hand side are permutable.
\eth

\Proof By Definition~\ref{def:comatrix} we have
\ben
\wh T(u)=\qdet T(u)\ts T^{-1}(u-n+1).
\end{equation*}
Taking the $nn$-th entry we come to
\ben
\qdet T(u)\ts \big(T^{-1}(u-n+1)\big)_{nn}=\wh t_{nn}(u).
\end{equation*}
Proposition~\ref{prop:mateco} gives
\ben
\qdet T(u)=\qdet T^{(n-1)}(u)\ts |T^{(n)}(u-n+1)|_{nn}.
\end{equation*}
Note that the factors here commute by the centrality of the quantum determinant.
An obvious induction completes the proof. \endproof

There is a generalization of this result providing a block factorization
of $\qdet T(u)$. For subsets ${\cal P}$ and ${\cal Q}$ in $\{1,\dots,n\}$
and an $n\times n$-matrix $X$
we shall denote by $X^{}_{{\cal P}{\cal Q}}$ the submatrix
of $X$ whose rows and columns are enumerated by 
${\cal P}$ and ${\cal Q}$ respectively.
Fix an integer $0\leq m\leq n$	and set ${\cal A}=\{1,\dots,m\}$
and ${\cal B}=\{m+1,\dots,n\}$. We let $T^*(u)$ denote the matrix $T^{-1}(-u)$.

\bth\label{thm:blocks}
We have the identity
\ben
\qdet T(u)\ts\qdet T^*(-u+n-1)_{{\cal A}{\cal A}}= \qdet T(u)_{{\cal B}{\cal B}}.
\end{equation*}
\eth

We keep the notation $\tpr_{ij}(u)$ for the matrix elements of the
matrix $T^{-1}(u)$.

\bpr\label{prop:matelinerse}
We have the relations
\ben
[t_{ij}(u),\tpr_{kl}(v)]=\frac{1}{u-v}
\Big(\delta_{kj}\sum_{a=1}^n t_{ia}(u)\ts\tpr_{al}(v)-
\delta_{il}\sum_{a=1}^n \tpr_{ka}(v)\ts t_{aj}(u)\Big).
\end{equation*}
In particular, the matrix elements of the matrices
$T(u)_{{\cal A}{\cal A}}$ and $T^*(v)_{{\cal B}{\cal B}}$
commute with each other.
\epr

\Proof It suffices to multiply the ternary relation \eqref{ternary}
by $T_2^{-1}(v)$ from both sides and equate the matrix elements. \endproof

\subsection{Quantum Sylvester theorem}\label{subsec:qst}

The following commutation relations between the quantum minors
generalize Proposition~\ref{prop:qminorel}.

\bpr\label{prop:qmr}
We have the relations
\ben
\begin{aligned}
{}&[{t\ts}^{a_1\cdots\ts a_k}_{b_1\cdots\ts b_k}(u),
{t\ts}^{c_1\cdots\ts c_l}_{d_1\cdots\ts d_l}(v)]=
\sum_{p=1}^{\min\{k,l\}}\frac{(-1)^{p-1}\ts p!}
{(u-v-k+1)\cdots (u-v-k+p)}
\\
{}&\sum_{\overset{\scl i_1<\cdots<i_p}
{\scl j_1<\cdots<j_p}}\left(
{t\ts}^{a_1\cdots\ts c_{j_1}\cdots\ts 
c_{j_p}\cdots\ts a_k}_{b_1\ \cdots\ \  b_k}(u)
{t\ts}^{c_1\cdots\ts a_{i_1}\cdots\ts 
a_{i_p}\cdots\ts c_l}_{d_1\ \cdots\ \  d_l}(v)
-{t\ts}^{c_1\ \cdots\ \  c_l}_{d_1\cdots\ts 
b_{i_1}\cdots\ts b_{i_p}\cdots\ts d_l}(v)
{t\ts}^{a_1\ \cdots\ \  a_k}_{b_1\cdots\ts 
d_{j_1}\cdots\ts d_{j_p}\cdots\ts b_k}(u)
\right).
\end{aligned}
\end{equation*}
Here the $p$-tuples
of upper indices $(a_{i_1},\dots, a_{i_p})$ and $(c_{j_1},\dots, c_{j_p})$
are respectively
interchanged in the first summand on the right hand
side while the $p$-tuples of lower indices
$(b_{i_1},\dots, b_{i_p})$ and
$(d_{j_1},\dots, d_{j_p})$
are interchanged in the second
summand.
\epr

As in the previous section, we fix an integer $m$ satisfying $1\leq m\leq n$.
For any indices $1\leq i,j\leq m$ introduce the following series
with coefficients in $\Y(n)$
\ben
{\wt t}_{ij}(u)={t\ts}^{i,m+1\cdots\ts n}_{j,m+1\cdots\ts n}(u)
\end{equation*}
and combine them into the matrix $\wt T(u)=(\wt t_{ij}(u))$.
The following is an analog of the classical Sylvester theorem.
As before, we denote ${\cal B}=\{m+1,\dots,n\}$.

\bth\label{thm:qs}
The mapping
\ben
t_{ij}(u)\mapsto {\wt t}_{ij}(u),\qquad 1\leq i,j\leq m,
\end{equation*}
defines an algebra homomorphism $\Y(m)\to\Y(n)$. Moreover,
one has the identity
\beql{qdettt}
\qdet {\wt T}(u)=\qdet T(u) \ts
\qdet T(u-1)^{}_{{\cal B}{\cal B}}\cdots\qdet
T(u-m+1)^{}_{{\cal B}{\cal B}}. 
\end{equation}
\eth

\Outline Using Proposition~\ref{prop:qmr} we check that the series ${\wt t}_{ij}(u)$
satisfy the Yangian defining relations \eqref{defrel} which proves the first claim.
The identity \eqref{qdettt} is derived by induction on $m$ from the relation
\ben
\wh T(u)^{}_{{\cal A}{\cal A}}\ts\wt T(u-m+1)=\qdet T(u)\ts\qdet
T(u-m+1)^{}_{{\cal B}{\cal B}},
\end{equation*}
where 
$\wh T(u)$ is the quantum comatrix; see Definition~\ref{def:comatrix}.
\endproof

\subsection{The centralizer construction}\label{subsec:cc}

Fix a nonnegative integer $m$ and
for any $n\geq m$ denote 
by $\g_m(n)$ the subalgebra in
$\gl_n$ spanned by the basis elements $E_{ij}$ with 
$m+1\leq i,j\leq n$. The subalgebra $\g_m(n)$ is isomorphic to
$\gl_{n-m}$.
Let $\Ar_m(n)$ denote the centralizer
of $\g_m(n)$ in the universal enveloping algebra $\Ar(n)=\U(\gl_n)$.
Let $\Ar(n)^0$ denote the centralizer of $E_{nn}$
in $\Ar(n)$ and let $\Ir(n)$ be the left ideal in $\Ar(n)$ generated by the
elements $E_{in}$, $i=1,\dots,n$. Then $\Ir(n)^0=\Ir(n)\cap \Ar(n)^0$
is a two-sided ideal in $\Ar(n)^0$ and one has a vector space
decomposition
\ben
\Ar(n)^0=\Ir(n)^0\oplus \Ar(n-1).
\end{equation*}
Therefore the projection of $\Ar(n)^0$ onto $\Ar(n-1)$ 
with the kernel $\Ir(n)^0$ is an algebra homomorphism. Its
restriction to the subalgebra $\Ar_m(n)$ defines
a filtration preserving homomorphism
\beql{proj}
\pi_n: \Ar_m(n)\to \Ar_m(n-1)
\end{equation}
so that one can define the algebra $\Ar_m$ as
the projective limit with respect to this sequence of homomorphisms
in the category of filtered algebras.

By the Harish-Chandra isomorphism \cite[Section~7.4]{d:ae}, the center
$\Ar_0(n)$ of  $\U(\gl_n)$ is naturally isomorphic to the algebra $\La^*(n)$
of polynomials in $n$ variables $\la_1,\dots,\la_n$ which are symmetric
in the shifted variables $\la_1,\la_2-1,\dots,\la_n-n+1$. So, in the case $m=0$
the homomorphisms $\pi_n$ are interpreted as the specialization
homomorphisms $\pi_n: \La^*(n)\to \La^*(n-1)$ such that
\beql{projm0}
\pi_n: f(\la_1,\dots,\la_n)\mapsto f(\la_1,\dots,\la_{n-1},0).
\end{equation}
The corresponding projective limit in the category of filtered
algebras is called the {\it algebra of shifted symmetric functions\/}
and denoted by $\La^*$. The elements of $\La^*$ are well-defined
functions on the set of all sequences $\la=(\la_1,\la_2,\dots)$ which
contain only a finite number of nonzero terms. 
The following families of elements of $\La^*$ are analogs of power
sums, elementary symmetric functions and complete symmetric functions:
\ben
\bal
p_m(\la)&=\sum_{k=1}^{\infty}\big((\la_k-k)^m-(-k)^m\big),\qquad m=1,2,\ldots,\\
1+\sum_{m=1}^{\infty}e_m(\la)t^m&=\prod_{k=1}^{\infty}
\frac{1+(\la_k-k)t}{1-kt},\\ 
1+\sum_{m=1}^{\infty}h_m(\la)t^m&=\prod_{k=1}^{\infty}
\frac{1+kt}{1-(\la_k-k)t}.
\eal
\end{equation*}
Each of the families $\{p_m\}$,
$\{e_m\}$, $\{h_m\}$ can be taken as a system of algebraically independent
generators of the algebra $\La^*$.
To summarize, we have the following.

\bpr\label{prop:shiftsym}
The algebra $\Ar_0$ is isomorphic to the algebra of shifted
symmetric functions $\La^*$.
\epr

Now consider the homomorphism \eqref{homtinv} and take its composition with an automorphism
of $\Y(n)$ given by \eqref{shift} to yield another homomorphism 
$\varphi_n : \Y(n)\to \Ar(n)$ such that
\ben
\varphi_n: T(u)\mapsto \Big(1-\frac{E}{u+n}\Big)^{-1}.
\end{equation*}
It follows from the defining relations \eqref{defrel} that the image of the
restriction of $\varphi_n$ to the subalgebra $\Y(m)$ is contained
in the centralizer $\Ar_m(n)$ thus yielding a homomorphism 
$\varphi_n : \Y(m)\to \Ar_m(n)$.

\bth\label{thm:cc}
For any fixed $m\geq 1$ the sequence
$(\varphi_n\ |\ n\geq m)$ defines an algebra embedding
$
\varphi: \Y(m)\hookrightarrow\Ar_m.
$
Moreover, one has
an isomorphism
\beql{isomcc}
\Ar_m=\Ar_0\ot \Y(m),
\end{equation}
where $\Y(m)$ is identified with its image under the embedding $\varphi$.
\eth

\Outline First we verify that the family  $(\varphi_n\ |\ n\geq m)$
is compatible with the chain of homomorphisms \eqref{proj}. 
Further, to prove the injectivity of $\varphi$ we consider
the corresponding commutative picture replacing $\Ar(n)$ with
its graded algebra $\gr\ts \Ar(n)\cong \Sr(\gl_n)$. This reduces
the task to the description of the invariants in the symmetric algebra
$\Sr(\gl_n)$ with respect to the action of the group $GL(n-m)$ corresponding
to the Lie algebra $\g_m(n)$. This can be done with the use of the classical invariant
theory \cite{w:cg} which also implies the decomposition \eqref{isomcc}. 
\endproof

A different embedding $\Y(m)\hookrightarrow\Ar_m$ can be constructed with the use
of the quantum Sylvester theorem; see Section~\ref{subsec:qst}.
Consider the homomorphism $\Y(m)\to\Y(n)$ provided by Theorem~\ref{thm:qs}
and take its composition with \eqref{eval}.
We obtain an algebra homomorphism $\psi_n:\Y(m)\to \Ar(n)$
given by
\beql{psin}
\psi_n: t_{ij}(u)\mapsto \qdet (1+Eu^{-1})^{}_{{\cal B}_i{\cal B}_j},
\end{equation}
where ${\cal B}_i$ denotes the set $\{i,m+1,\dots,n\}$.
Corollary~\ref{cor:qmincent} implies that this image commutes
with the elements of the subalgebra $\g_m(n)$ and so,
\eqref{psin} defines a homomorphism $\psi_n:\Y(m)\to \Ar_m(n)$.
Furthermore, the family of homomorphisms $(\psi_n\ |\ n\geq m)$
is obviously compatible with the projections \eqref{proj}
and thus defines an algebra homomorphism $\psi: \Y(m)\to \Ar_m$.
Denote by ${\wt \Ar}_0$ the projective limit of the sequence
of the centers of the universal enveloping algebras
$\U(\g_m(n))$, where $n=m,m+1,\dots$, defined by the corresponding
homomorphisms \eqref{projm0}. By Proposition~\ref{prop:shiftsym},
${\wt \Ar}_0$ is isomorphic to the algebra of shifted symmetric
functions in the variables $\la_{m+1},\la_{m+2},\dots$.
The following is an analog of Theorem~\ref{thm:cc}.

\bth\label{thm:mcc}
The homomorphism $\psi: \Y(m)\to \Ar_m$ is injective. Moreover, one has
an isomorphism
\ben
\Ar_m={\wt \Ar}_0\ot \Y(m),
\end{equation*}
where $\Y(m)$ is identified with its image under the embedding $\psi$.
\eth

\subsection{Commutative subalgebras}\label{subsec:cs}

Here we use the notation introduced in Section~\ref{subsec:matrix}.
Consider the algebra \eqref{multitp} with $m=n$.
Fix an $n\times n$ matrix $C$ with entries in $\CC$
and for any $1\leq k\leq n$
introduce the series $\tau_k(u,C)$ with coefficients in $\Y(n)$ by
\ben
\tau_k(u,C)=\tr A_n\ts {T}_1(u)\cdots {T}_k(u-k+1)\ts C_{k+1}\cdots C_n,
\end{equation*}
where $A_n$ is the antisymmetrizer defined by \eqref{antisym} and the trace
is taken over all $n$ copies of $\End\CC^n$. 

\bth\label{thm:commsub}
All the coefficients of the series $\tau^{}_1(u,C),\dots,\tau^{}_n(u,C)$ commute
with each other.
Moreover, if the matrix $C$ has simple spectrum then 
the coefficients at $u^{-1},u^{-2},\dots$ of these series
are algebraically independent and generate a 
maximal commutative subalgebra of the Yangian $\Y(n)$.
\eth

Consider the epimorphism $\pi:\Y(n)\to \U(\gl_n)$ defined in \eqref{eval}.
Clearly, the coefficients of the images
of the series $\tau_k(u,C)$, $k=1,\dots,n$ under $\pi$ form a commutative subalgebra
${\mathcal C}\subseteq \U(\gl_n)$. 

\bth\label{thm:commsubu}
If the matrix $C$ has simple spectrum then the subalgebra
${\mathcal C}$ of \ts $\U(\gl_n)$ is maximal commutative.
\eth

\subsection*{Bibliographical notes}\label{subsec:bibgln}

{\bf \ref{subsec:matrix}.}
For the origins of the $RTT$ relation and associated 
$R$-matrix formalism see for instance the papers Takhtajan--Faddeev~\cite{tf:qi},
Kulish--Sklyanin~\cite{ks:qs}, Reshetikhin--Takhtajan--Faddeev~\cite{rtf:ql}.
The statistical mechanics background of quantum groups is also explained in
the book by Chari and Pressley~\cite[Chapter~7]{cp:gq}.

\noindent
{\bf \ref{subsec:pbw}.} The Poincar\'e--Birkhoff--Witt
theorem for general Yangians is due to Drinfeld
(unpublished). Another proof was given by Levendorski\u\i~\cite{l:pb}.
The details of the proof outlined here
can be found in \cite{mno:yc}. It follows the approach of 
Olshanski\u\i~\cite{o:ri}.

\noindent{\bf \ref{subsec:qdet}.}
The definition of the quantum determinant $\qdet T(u)$ (in the case $n=2$)
originally appeared in Izergin--Korepin~\cite{ik:lm} .
The basic ideas and formulas associated with the
quantum determinant for an arbitrary $n$
are contained in Kulish--Sklyanin's survey paper
\cite{ks:qs}. Detailed proofs are given in \cite{mno:yc}.
Proposition~\ref{prop:coprqm} is contained e.g. in 
Iohara~\cite{i:br} and Nazarov--Tarasov~\cite{nt:ry}.

\noindent{\bf \ref{subsec:sln}.}
By the general approach of Drinfeld~\cite{d:qg}, 
the Yangian for $\sll_n$ should be
defined as a quotient algebra of $\Y(n)$. The fact that it can also be
realized as a (Hopf)
subalgebra of $\Y(n)$ was observed by Olshanski~\cite{mno:yc}.

\noindent{\bf \ref{subsec:3real}.}
For any simple Lie algebra $\agot$ the Yangian $\Y(\agot)$
was defined by Drinfeld~\cite{d:ha, d:nr}.
The two definitions given here are particular cases
for $\agot=\sll_2$. 
The formulas for the coproduct and antipode
in Definition~\ref{def:ydefgf2} are due to the author;
see e.g. Khoroshkin--Tolstoy~\cite{kt:yd}. These formulas were employed
in \cite{kt:yd} in the construction of the double of the Yangian.
These results were generalized to the Yangian $\Y(\sll_3)$
by Soloviev~\cite{s:cw}, and to $\Y(\sll_n)$ by Iohara~\cite{i:br}.

\noindent{\bf \ref{subsec:qlf}.} 
The series $z(u)$ was introduced by Nazarov~\cite{n:qb}.
The quantum Liouville formula is also due to him.
The argument given in
Section~\ref{subsec:qlf} is a 
simplified version of his $R$-matrix proof \cite{mno:yc}.

\noindent{\bf \ref{subsec:qdf}.}
A general theory of quasi-determinants of matrices
over noncommutative rings is developed by 
Gelfand and Retakh~\cite{gr:dm, gr:tn}.
Various analogs of the classical theorems for such determinants
are given. 
Quasi-determinant factorizations of quantum determinants
for the quantized algebra $GL_q(n)$
are obtained in those papers; see also
Krob and Leclerc~\cite{kl:mi}.
In a more general context of Hopf algebras such factorizations
are constructed by Etingof and Retakh~\cite{er:qd}.

\noindent{\bf \ref{subsec:qst}.}
An analog of Sylvester's theorem for
the algebra $GL_q(n)$
was given
by Krob and Leclerc~\cite{kl:mi} with the use 
of the quasi-determinant version
of this theorem due to Gelfand and Retakh~\cite{gr:dm}.
The approach of \cite{kl:mi} works for 
the Yangians as well. The proof outlined here follows 
\cite{m:yt} where a proof of
Proposition~\ref{prop:qmr} is given. The latter result and some other
quantum minor relations are known 
to specialists as `folklore theorems'.
Some more quantum analogs of the classical minor relations
are collected in Iohara~\cite{i:br}. 

\noindent{\bf \ref{subsec:cc}.} 
The	centralizer construction is due to Olshanski~\cite{o:ea, o:ri}.
The modified version (Theorem~\ref{thm:mcc}) based on the quantum Sylvester theorem
is given in \cite{m:yt}. The algebra $\La^*$ of shifted symmetric functions
is studied in detail by Okounkov and Olshanski~\cite{oo:ss}.

\noindent{\bf \ref{subsec:cs}.} The
commutative subalgebras in the Yangian originate from
the integrable models in statistical mechanics (specifically, 
from the {\it six vertex\/}
or {\it XXX\/} model); see Baxter~\cite{b:es}. 
The common eigenvectors of such a commutative subalgebra
in certain standard Yangian modules can be constructed
by a special procedure called the {\it algebraic Bethe ansatz\/}; 
see Faddeev~\cite{f:im},
Kulish--Sklyanin~\cite{ks:qs}, Kulish--Reshetikhin~\cite{kr:dg},
Kirillov--Reshetikhin~\cite{kr:yb, kr:ba}.
Theorems~\ref{thm:commsub} and \ref{thm:commsubu} are proved by 
Nazarov and Olshanski~\cite{no:bs}.

\medskip
Finite-dimensional irreducible representations of the Yangians were classified
by Drinfeld~\cite{d:nr} with the use of the particular case of $\Y(\gl_2)$
considered earlier by Tarasov~\cite{t:sq, t:im}. A detailed exposition of these
results for $\Y(\gl_n)$ is contained in \cite{m:fd}.
Cherednik~\cite{c:ni,c:qg} used the Yangians to `materialize' the second Weyl
character formula. Nazarov~\cite{n:qb, n:yc} and Okounkov~\cite{o:qi}
employed the Yangian techniques to obtain remarkable immanant analogs
of the classical Capelli identity \cite{c:zc, c:ot}.
An explicit construction of all representations of $\Y(\gl_2)$ as tensor products of
the evaluation modules is given by Tarasov~\cite{t:im} and
Chari--Pressley~\cite{cp:yr, cp:yt}. In particular, this provides
an irreducibility criterion of tensor products
of the $\Y(\gl_2)$ evaluation modules. 
A generalization of this criterion
to $\Y(\gl_n)$ with an arbitrary $n$ was given in \cite{m:ic}; see also
Leclerc--Nazarov--Thibon~\cite{lnt:ir}. 
An important part of the criterion is
the {\it binary property\/} established by
Nazarov and Tarasov~\cite{nt:it}; see also 
Kitanine, Maillet and Terras~\cite{kmt:ff, mt:qi}.
A character formula for an arbitrary finite-dimensional irreducible 
representation of $\Y(\gl_n)$ is given by Arakawa~\cite{a:df}
with the use of the Drinfeld functor; see Drinfeld~\cite{d:da}.
This formula is also implied by the earlier results of Ginzburg--Vasserot~\cite{gv:lr}
combined with the work of Lusztig; see Nakajima~\cite{n:qv}, Varagnolo~\cite{v:qv}.
The irreducible characters are expressed in terms of those for the 
`standard tensor product modules' via
the Kazhdan--Lusztig polynomials. Bases of Gelfand--Tsetlin type
for `generic' representations of $\Y(\gl_n)$ are constructed in \cite{m:gt}.
More general class of `tame' Yangian modules was introduced and explicitly
constructed by Nazarov and Tarasov~\cite{nt:ry}.
The earlier works of Nazarov and Tarasov~\cite{nt:yg} and the author~\cite{m:gt}
provide different constructions of the well-known Gelfand--Tsetlin basis
for representations of $\gl_n$. A surprising connection
of the Yangian $\Y(\gl_n)$ with the finite $\mathcal W$-algebras was discovered
by Ragoucy and Sorba~\cite{rs:yf, rs:yr}; see also Briot and Ragoucy~\cite{br:rp}.
The Yangian actions on certain modules over the affine Lie algebras 
were constructed by Uglov~\cite{u:ss}.

\section{The twisted Yangians}\label{sec:twyang}
\setcounter{equation}{0}

Here we describe the structure of the twisted Yangians corresponding to
the orthogonal and symplectic Lie algebras $\oa^{}_N$ and $\spa^{}_N$.
We consider both cases simultaneously, unless otherwise stated.

\subsection{Defining relations}\label{subsec:twdefrel}

Given a positive integer $N$, we number the rows and columns
of $N\times N$ matrices by the indices
$\{-n,\dots,-1,0,1,\dots n\}$ if $N=2n+1$, and by
$\{-n,\dots,-1,1,\dots n\}$ if $N=2n$.
Similarly, in the latter case the range of
indices $-n\leq i,j\leq n$ will exclude $0$.
It will be convenient to use
the symbol $\theta_{ij}$ which is defined by
\ben
\theta_{ij}=\begin{cases} 1\quad&\text{in the orthogonal case},\\
\sgn i\cdot \sgn j\quad&\text{in the symplectic case}.\end{cases}
\end{equation*}
Whenever the double sign $\pm{}$ or $\mp{}$ occurs,
the upper sign corresponds to the orthogonal case and the lower sign to
the symplectic case.
By $X\mapsto X^t$ we will denote the matrix 
transposition such that	for the matrix units we have
\beql{transptw}
(e_{ij})^t=\theta_{ij}e_{-j,-i}.
\end{equation}
Introduce the following elements of the Lie algebra $\gl^{}_N$:
\ben
F_{ij}=E_{ij}-\theta_{ij}E_{-j,-i},\qquad -n\leq i,j\leq n.
\end{equation*}
The Lie subalgebra of $\gl^{}_N$ spanned 
by the elements $F_{ij}$ is a realization of a simple Lie algebra $\g^{}_n$
of rank $n$ (see e.g. \cite{h:il}).	In the orthogonal case
$\g^{}_n$ is of type  $D_n$ or $B_n$
if $N=2n$ or $N=2n+1$, respectively.
In the symplectic case, $N=2n$ and $\g^{}_n$ is of type $C_n$.
Thus,
\ben
\g^{}_n=\oa^{}_{2n},\quad  \oa^{}_{2n+1}\quad\text{or}\quad \spa^{}_{2n}.
\end{equation*}

\bde\label{def:twyang}
Each of the {\em twisted Yangians\/} $\Y^{+}(2n)$, $\Y^{+}(2n+1)$ and $\Y^{-}(2n)$
corresponding to the Lie algebras $\oa^{}_{2n}$, $\oa^{}_{2n+1}$ and
$\spa^{}_{2n}$, respectively, is a unital associative algebra
with generators $s_{ij}^{(1)},\ s_{ij}^{(2)},\dots$ where $-n\leq i,j\leq n$,
and the defining relations written in terms of
the generating series
\ben
s_{ij}(u)=\delta_{ij}+ s_{ij}^{(1)}u^{-1}+s_{ij}^{(2)}u^{-2}+\cdots
\end{equation*}
as follows
\ben
\bal
(u^2-v^2)\ts
[s_{ij}(u),s_{kl}(v)]=(u+v)\ts&\big(s_{kj}(u)s_{il}(v)-s_{kj}(v)s_{il}(u)\big)\\
{}-(u-v)\ts&\big(\theta_{k,-j}s_{i,-k}(u)s_{-j,l}(v)-
\theta_{i,-l}s_{k,-i}(v)s_{-l,j}(u)\big)\\
{}+\theta_{i,-j}&\big(s_{k,-i}(u)s_{-j,l}(v)-s_{k,-i}(v)s_{-j,l}(u)\big)
\eal
\end{equation*}
and
\beql{symme}
\theta_{ij}s_{-j,-i}(-u)=s_{ij}(u)\pm \frac{s_{ij}(u)-s_{ij}(-u)}{2u}. 
\end{equation}
\ede

These relations	can be conveniently presented in an equivalent matrix form
analogous to \eqref{ternary}. For this introduce the transposed $R^{\tss t}(u)$
for the Yang $R$-matrix \eqref{rmatrix} (correcting the indices 
for the matrix elements appropriately) by
\ben
R^{\tss t}(u)=1-Q\ts u^{-1}, \qquad Q=\sum_{i,j=-n}^n e^t_{ij}\ot e_{ji}.
\end{equation*}
Furthermore, denote by $S(u)$ the $N\times N$ matrix whose 
$ij$-th entry is $s_{ij}(u)$. As in \eqref{tmatrix} we regard
$S(u)$ as an element of the algebra $\Y^{\pm}(N)\ot\End\CC^N$ given by
\ben
S(u)=\sum_{i,j=-n}^n s_{ij}(u)\ot e_{ij}.
\end{equation*}

\bpr\label{prop:quatern}
The defining relations for the twisted Yangian $\Y^{\pm}(N)$ are equivalent to
the quaternary relation
\beql{quatern}
R(u-v)S_1(u)R^{\tss t}(-u-v)S_2(v)=S_2(v)R^{\tss t}(-u-v)S_1(u)R(u-v)
\end{equation}
and the symmetry relation
\beql{symmet}
S^t(-u)=S(u)\pm \frac{S(u)-S(-u)}{2u}.
\end{equation}
\epr

The following relation between the twisted Yangians and the corresponding 
classical Lie algebras plays a key role in many applications; 
cf. Proposition~\ref{prop:eval}.

\bpr\label{prop:evaltw}
The mapping
\beql{evaltw}
\pi : s_{ij}(u)\mapsto \delta_{ij}+F_{ij}\Big(u\pm\frac12\Big)^{-1}
\end{equation}
defines an algebra epimorphism $\Y^{\pm}(N)\to\U(\g^{}_n)$. Moreover,
\ben
F_{ij}\mapsto s_{ij}^{(1)}
\end{equation*}
is an embedding $\U(\g^{}_n)\hookrightarrow \Y^{\pm}(N)$.
\epr

\subsection{Embedding into the Yangian}\label{subsec:embed}

We keep the notation $t_{ij}^{(r)}$ for the generators of the Yangian $\Y(N)$.
However, in accordance with the above, we now let the indices $i,j$ run over
the set $\{-n,\dots,n\}$.  Also, the matrix transposition is now understood in the
sense \eqref{transptw}.

\bth\label{thm:embed}
The mapping
\beql{embedtwy}
S(u)\mapsto T(u)\ts T^{\tss t}(-u)
\end{equation}
defines an embedding $\Y^{\pm}(N)\hookrightarrow\Y(N)$.
\eth

\Outline It is straightforward to verify that the matrix $T(u)\ts T^{\tss t}(-u)$
satisfies both relations \eqref{quatern} and \eqref{symmet}.
To show that the homomorphism \eqref{embedtwy} is injective we
use the corresponding homomorphism of the graded algebras
$\gr_1 \Y^{\pm}(N)\to \gr_1 \Y(N)$, where $\gr_1 \Y^{\pm}(N)$ is defined
by setting $\deg_1 s_{ij}^{(r)}=1$. Then apply Corollary~\ref{cor:polalg}.
\endproof

This result allows us to regard the twisted Yangian as a subalgebra
of $\Y(N)$. The following is an analog of the 
Poincar\'e--Birkhoff--Witt theorem for the algebra $\Y^{\pm}(N)$.

\bco\label{cor:pbwtw}
Given an arbitrary linear order on the set of
the generators
\ben
s_{ij}^{(2k)},\quad i+j\leq 0;\qquad  s_{ij}^{(2k-1)},\quad
i+j<0;\qquad k=1,2,\dots,
\end{equation*}
in the case of $\Y^{+}(N)$, and the set of the generators
\ben
s_{ij}^{(2k)},\quad i+j< 0;\qquad  s_{ij}^{(2k-1)},\quad
i+j\leq 0;\qquad k=1,2,\dots,
\end{equation*}
in the case of $\Y^{-}(N)$,
any element of the algebra $\Y^{\pm}(N)$ is uniquely written
as a linear combination of ordered monomials in the generators.
\eco

\bpr\label{prop:coideal}
The subalgebra $\Y^{\pm}(N)$ is a left coideal of $\Y(N)$,
\ben
\Delta\big(\Y^{\pm}(N)\big)\subseteq \Y(N)\ot \Y^{\pm}(N)
\end{equation*}
\epr

\Proof This follows from the explicit formula
\ben
\Delta: s_{ij}(u)\mapsto \sum_{a,b=-n}^n 
\theta_{bj}t_{ia}(u)t_{-j,-b}(-u)\ot s_{ab}(u).
\end{equation*}

The restrictions of some of the automorphisms and anti-automorphisms
of $\Y(N)$ described in Section~\ref{subsec:auto} preserve the subalgebra
$\Y^{\pm}(N)$ and so we have the following.

\bpr\label{prop:autotw}
The mapping 
\ben
S(u)\mapsto S^{\tss t}(u)
\end{equation*}
defines an anti-automorphism of $\Y^{\pm}(N)$.
For any formal series $g(u)\in 1+u^{-2}\CC[[u^{-2}]]$ the mapping
\beql{mulsuu}
S(u)\mapsto g(u)\ts S(u)
\end{equation}
defines an automorphism of $\Y^{\pm}(N)$.
\epr

\subsection{Sklyanin determinant}\label{subsec:sklyanin}

Here we use the notation of Section~\ref{subsec:qdet} with the usual convention
on the matrix element indices. Define $R(u_1,\dots,u_m)$ by \eqref{prodrma}
and set
\ben
S_i=S_i(u_i),\quad 1\leq i\leq m \qquad \text{and}\qquad
R^{\tss t}_{ij}=R^{\tss t}_{ji}=R^{\tss t}_{ij}(-u_i-u_j),\quad 1\leq i<j \leq m,
\end{equation*}
where $S_a(u)$ and $R^{\tss t}_{ab}(u)$ are defined by an obvious analogy
with \eqref{tku} and \eqref{ars}.
For an arbitrary permutation $(p_1,\ldots,p_m)$ of the numbers
$1,\ldots,m$, we abbreviate
\ben
\langle S_{p_1},\ldots,S_{p_m}\rangle=S_{p_1}(R^{\tss t}_{p_1p_2}\cdots
R^{\tss t}_{p_1p_m}) S_{p_2}(R^{\tss t}_{p_2p_3}\cdots R^{\tss t}_{p_2p_m})\cdots S_{p_m}.
\end{equation*}

\bpr\label{prop:fundamtw} We have the identity
\ben
R(u_1,\ldots,u_m)\langle S_1,\ldots,S_m\rangle=\langle
S_m,\ldots,S_1\rangle R(u_1,\ldots,u_m).
\end{equation*}
\epr

Now take $m=N$ and specialize the variables $u_i$ as
\ben
u_i=u-i+1,\qquad i=1,\dots,N.
\end{equation*}
By Propositions~\ref{prop:rmanti} and \ref{prop:fundamtw} we have
\beql{antw}
A_N	\ts\langle S_1,\ldots,S_N\rangle=\langle
S_N,\ldots,S_1\rangle\ts A_N.
\end{equation}

\bde\label{def:sdet}
The {\ts\em Sklyanin determinant\/} is the formal series
\ben
\sdet S(u)=1+c_1 u^{-1}+c_2 u^{-2}+\cdots\in\Y^{\pm}(N)[[u^{-1}]]
\end{equation*}
such that the element \eqref{antw} equals $A_N\ts\sdet S(u)$.
\ede

In the next theorem we regard $\Y^{\pm}(N)$ as a subalgebra of $\Y(N)$.

\bth\label{thm:sdetqdet} We have
\ben
\sdet S(u)=\gamma^{}_n(u)\ts\qdet T(u)\ts\qdet T(-u+N-1),
\end{equation*}
where
\ben
\gamma^{}_n(u)=\begin{cases}\  1\qquad&\text{for}\quad \Y^+(N),\\
             \dfrac{2u+1}{2u-2n+1}\qquad&\text{for}\quad \Y^-(2n).
			   \end{cases}
\end{equation*}
\eth

There is an explicit formula for the Sklyanin determinant in terms of
the generators $s_{ij}(u)$. It uses a special map of the symmetric groups
\beql{mapsym}
\pi^{}_N:\Sym_N\to \Sym_{N},\qquad p\mapsto p'
\end{equation}
defined by an inductive procedure. First of all, $p'(N)=N$ so that
$p'$ can be regarded as an element of $\Sym_{N-1}$.
Given a set of positive integers 
$\om_1<\cdots<\om_N$ we 
regard $\Sym_N$ as the group of their permutations.
If $N=2$ we define $\pi^{}_2$ as the only map $\Sym_2\to \Sym_{1}$.
For $N>2$ define a map from the set of ordered pairs $(\om_k,\om_l)$
with $k\ne l$ into itself by the rule
\beql{ordpair}
\begin{alignedat}{2}
(\om_k,\om_l)&\mapsto (\om_l,\om_k),&&\qquad k,l<N,\\
(\om_k,\om_N)&\mapsto (\om_{N-1},\om_k),&&\qquad k<N-1,\\
(\om_N,\om_k)&\mapsto (\om_k,\om_{N-1}),&&\qquad k<N-1,\\
(\om_{N-1},\om_N)&\mapsto (\om_{N-1},\om_{N-2}),\\
(\om_{N},\om_{N-1})&\mapsto (\om_{N-1},\om_{N-2}).
\end{alignedat}
\end{equation}
Let $p=(p^{}_1,\dots,p^{}_N)$ be a permutation of the indices
$\om_1,\dots,\om_N$. Define its image $p'=(q^{}_1,\dots,q^{}_{N-1})$ under
the map $\pi^{}_N$ as follows. First take $(q^{}_1,q^{}_{N-1})$ as the image
of the ordered pair $(p^{}_1,p^{}_N)$ under the map \eqref{ordpair}.
Then define $(q^{}_2,\dots,q^{}_{N-2})$ as the image of $(p^{}_2,\dots,p^{}_{N-1})$
under $\pi^{}_{N-2}$ where $\Sym_{N-2}$ is regarded as the group 
of permutations of the family of indices obtained from $(\om_1,\dots,\om_N)$
by deleting $p^{}_1$ and $p^{}_N$. To describe the combinatorial properties
of the map $\pi^{}_N$ introduce the signless Stirling numbers of the first kind
$c(N,k)$ by
\ben
\sum_{k=1}^N c(N,k)\ts x^k=x(x+1)\cdots (x+N-1).
\end{equation*}
We also need the Boolean posets $B_n$. The elements of $B_n$
are subsets of $\{1,\dots,n\}$ with the usual set inclusion 
as the partial ordering. Equip the symmetric group $\Sym_N$
with  the standard Bruhat order.

\bth\label{thm:stirling}
Each fiber of the map $\pi^{}_N:\Sym_N\to\Sym_{N}$ is an interval in $\Sym_N$
with respect to the Bruhat order, isomorphic to the Boolean poset $B_k$ for some $k$.
Moreover, for any $k$ the number of intervals isomorphic to $B_k$
is the Stirling number $c(N-1,k)$.
\eth

Below we denote the matrix elements of the transposed matrix $S^{\tss t}(u)$
by $s^t_{ij}(u)$. For any permutation $p\in \Sym_N$ we denote by $p'$
its image under the map $\pi^{}_N$.

\bth\label{thm:sklproj}
Let $(a^{}_1,\dots,a^{}_N)$ be an arbitrary permutation of the set of  indices
$(-n,-n+1,\dots,n)$. Then
\ben
\bal
\sdet S(u)=(-1)^n\gamma^{}_n(u)\sum_{p\in\Sym_N}
&\sgn p\tss p'\cdot s^t_{-a^{}_{p(1)},\ts a^{}_{p'(1)}}(-u)\cdots 
s^t_{-a^{}_{p(n)},\ts a^{}_{p'(n)}}(-u+n-1)\\
{}&{}\times s^{}_{-a^{}_{p(n+1)},\ts a^{}_{p'(n+1)}}(u-n)\cdots
s^{}_{-a^{}_{p(N)},\ts a^{}_{p'(N)}}(u-N+1)
\eal
\end{equation*}
and also
\ben
\bal
{}=(-1)^n\gamma^{}_n(u)\sum_{p\in\Sym_N}
&\sgn p\tss p'\cdot s^{}_{-a^{}_{p'(1)},\ts a^{}_{p(1)}}(u-N+1)\cdots
s^{}_{-a^{}_{p'(n)},\ts a^{}_{p(n)}}(u-N+n)\\
{}&{}\times s^t_{-a^{}_{p'(n+1)},\ts a^{}_{p(n+1)}}(-u+N-n-1)\cdots
s^t_{-a^{}_{p'(N)},\ts a^{}_{p(N)}}(-u).
\eal
\end{equation*}
\eth

\Outline The key idea is to use the symmetry relation \eqref{symmet}
in order to
eliminate the intermediate factors
$R^{\tss t}_{ij}$ in the expression \eqref{antw} which defines 
the Sklyanin determinant.
\endproof

\bex\label{ex:N=2}
In the case $N=2$ we have
\ben
\bal
\sdet S(u)=
\frac{2u+1}{2u\pm 1}&\big(s^{}_{-1,-1}(u-1)s^{}_{-1,-1}(-u)\mp
s^{}_{-1,1}(u-1)s^{}_{1,-1}(-u)\big)\\
{}=\frac{2u+1}{2u\pm 1}&\big(s^{}_{1,1}(-u)s^{}_{1,1}(u-1)\mp 
s^{}_{1,-1}(-u)s^{}_{-1,1}(u-1)\big).
\eal
\end{equation*}
\eex

\subsection{The center of the twisted Yangian}\label{subsec:centertw}

Theorem~\ref{thm:sdetqdet} implies the following relation between the coefficients
of the Sklyanin determinant
\ben
\gamma^{}_n(u)\ts \sdet S(-u+N-1)=\gamma^{}_n(-u+N-1)\ts \sdet S(u).
\end{equation*}
In particular, the odd coefficients of $\sdet S(u)$ can be expressed in terms
of the even.

\bth\label{thm:centertw}
All the coefficients of the series $\sdet S(u)$
belong to the center of the algebra $\Y^{\pm}(N)$. 
Moreover, the even coefficients $c_2,c_4,\dots$ are algebraically
independent and generate the center of $\Y^{\pm}(N)$.
\eth

\Outline The first assertion is immediate from Theorems~\ref{thm:qdcenter}
and \eqref{thm:sdetqdet}. Alternatively, one can also prove the centrality
of $\sdet S(u)$ by using the matrix form of the defining relations
of $\Y^{\pm}(N)$. To prove the second assertion introduce a filtration
on $\Y^{\pm}(N)$ by setting $\deg_2 s_{ij}^{(r)}=r-1$.
To describe the corresponding graded algebra $\gr_2 \Y^{\pm}(N)$
consider the involution	$\sigma$ on the Lie algebra $\gl^{}_N$
defined by $\sigma:A\mapsto -A^t$. Then the fixed point subalgebra
is the classical Lie algebra $\g_n$. Introduce the twisted polynomial
current Lie algebra $\gl^{}_N[x]^{\sigma}$ by \eqref{twpolc}.
Then we have an algebra isomorphism
\ben
\gr_2 \Y^{\pm}(N)\cong \U(\gl^{}_N[x]^{\sigma}).
\end{equation*}
Identifying these algebras we find that
the image of the coefficient $c_{2m}$ in the $(2m-1)$-st component of
$\gr_2 \Y^{\pm}(N)$ coincides with $Z\ts x^{2m-1}$ where
\ben
Z=E_{-n,-n}+E_{-n+1,-n+1}+\cdots+E_{n,n}.
\end{equation*}
To complete the proof of the theorem we use the fact that
the center of the algebra $\U(\gl^{}_N[x]^{\sigma})$ is generated
by the elements $Z\ts x^{2m-1}$ with $m\geq 1$.	\endproof

\subsection{The special twisted Yangian}\label{subsec:spectw}

In the next definition we regard $\Y^{\pm}(N)$ as a 
subalgebra of the Yangian $\Y(N)$. Recall the definition of the Yangian
$\Y(\sll_N)$ from Section~\ref{subsec:sln}.

\bde\label{def:spectw}
The {\em special twisted Yangian\/} $\SY^{\pm}(N)$ is the subalgebra
of $\Y^{\pm}(N)$ defined by
\ben
\SY^{\pm}(N)=\Y(\sll^{}_N)\cap \Y^{\pm}(N).
\end{equation*}
\ede

Equivalently, $\SY^{\pm}(N)$ is the subalgebra of $\Y^{\pm}(N)$
which consists of the elements stable under all automorphisms
of the form \eqref{mulsuu}.	The following result is implied by
Theorem~\ref{thm:decoyn}.

\bth\label{thm:twdeco}
The algebra  $\Y^{\pm}(N)$ is isomorphic
to
the tensor product of its center $\Z^{\pm}(N)$ and the subalgebra
$\SY^{\pm}(N)$,
\ben
\Y^{\pm}(N)=\Z^{\pm}(N)\otimes \SY^{\pm}(N).
\end{equation*}
In particular, the center of $\SY^{\pm}(N)$ is trivial.
\eth

\bco\label{cor:coidts}
The subalgebra $\SY^{\pm}(N)$ of $\ts\Y(\sll^{}_N)$ is a left coideal.
\eco

\subsection{The quantum Liouville formula}\label{subsec:liouvtw}

Define the series $\zeta (u)$ with coefficients in $\Y^{\pm}(N)$ by
\ben
\zeta(u)^{-1}=\frac 1N\ts\tr
\Big\{\Big(\frac{2u-N}{2u-N\pm1}S^{\tss t}(-u)
\pm\frac1{2u-N\pm1}S(-u)\Big)S^{-1}(u-N)\Big\}.
\end{equation*}

\bth\label{thm:liouvtw}
We have the relation
\ben
\zeta(u)=\varepsilon^{}_n(u)\ts\frac{\sdet S(u-1)}{\sdet S(u)},
\end{equation*}
where $\varepsilon^{}_n(u)=\gamma^{}_n(u)\ts\gamma^{}_n(u-1)^{-1}$.
\eth

\Outline The quaternary relation \eqref{quatern} implies
\ben
Q\ts S_1^{-1}(-u)\ts R(2u-N)\ts S_2(u-N)=S_2(u-N)\ts R(2u-N)\ts S_1^{-1}(-u)\ts Q.
\end{equation*}
It is deduced from the definition of $\zeta(u)$ that this expression
coincides with $\zeta(u)\ts Q$ up to a scalar function. Combining this with
the matrix definition \eqref{qresi} of $z(u)$ we obtain
\ben
\zeta(u)=z(u)\ts z(-u+N)^{-1}.
\end{equation*}
Now the claim follows from Theorems~\ref{thm:liouville} and \ref{thm:sdetqdet}.
\endproof

\bco\label{cor:zetagen}
The coefficients of the series $\zeta(u)$ generate the center of $\Y^{\pm}(N)$.
\eco

\subsection{Factorization of the Sklyanin determinant}
\label{subsec:sdf}

Let $1\leq m\leq n$.
Denote by $S^{(m)}(u)$ the submatrix of $S(u)$
corresponding the rows and columns enumerated by $-m,-m+1,\dots,m$,
and by $\wt S^{(m)}(u)$ the submatrix of $S^{(m)}(u)$ obtained by removing
the row and column enumerated by $-m$. Set
\ben
c(u)=\frac{1}{\gamma^{}_n(u+N/2-1/2)}\ts
\sdet S(u+N/2-1/2).
\end{equation*}
Then by Theorem~\ref{thm:sdetqdet}, $c(u)$ is an even formal series in $u^{-1}$,
with coefficients in the center of the $\Y^{\pm}(N)$.
We shall use the quasi-determinants introduced in \eqref{def:quasidet}.

\bth\label{thm:sdetdec} If $N=2n$ then
\ben
\bal
c(u)=|\wt S^{(1)}(-u-1/2)|_{11}\cdot&\ts|S^{(1)}(u-1/2)|_{11}\\
\cdots &\ts|\wt S^{(n)}(-u-n+1/2)|_{nn}\cdot|S^{(n)}(u-n+1/2)|_{nn}.
\eal
\end{equation*}
If $N=2n+1$ then
\ben
\bal
c(u)=s_{00}(u)\cdot|\wt S^{(1)}(-u-1)|_{11}\cdot|S^{(1)}(u-1)|_{11}
\cdots |\wt S^{(n)}(-u-n)|_{nn}\cdot|S^{(n)}(u-n)|_{nn}.
\eal
\end{equation*}
Moreover, the factors on the right side of each
expression are permutable.
\eth

\Outline Define the {\it Sklyanin comatrix\/} $\wh S(u)=
(\wh s_{ij}(u))$ by the formula
\beql{sklco}
\wh S(u)\ts S(u-N+1)=\sdet S(u).
\end{equation}
Taking the $nn$-th entry gives
\ben
\sdet S(u)=\wh s_{nn}(u)\ts |S(u-N+1)|_{nn}.
\end{equation*}
Then proceed by induction with the use of the formula
\ben
\wh s_{nn}(u)=\frac{2u+1}{2u\pm 1}\ts
|\wt S^{(n)}(-u)|_{nn}\ts\sdet S^{(n-1)}(u-1),
\end{equation*}
which is deduced from the definition of $\sdet S(u)$.\endproof

\bpr\label{prop:sklauto}
The mapping
\ben
S(u)\mapsto \gamma_N(u)\ts \wh S(-u+\frac N2-1),
\end{equation*}
defines an automorphism of the algebra $\Y^{\pm}(N)$.
\epr

Consider the algebra $\wt{\Y}^{\pm}(N)$ which is defined exactly as the twisted
Yangian  $\Y^{\pm}(N)$ (Definition~\ref{def:twyang}) 
but with the symmetry relation \eqref{symme} dropped. We denote
the generators of $\wt{\Y}^{\pm}(N)$ by the same symbols $s_{ij}^{(r)}$.
Note that the definition of the Sklyanin determinant
does not use the symmetry relation; see Definition~\ref{def:sdet}.
Therefore, we can define $\sdet S(u)$ in the same way for
$\wt{\Y}^{\pm}(N)$. One can show that all coefficients of $\sdet S(u)$
belong to the center of this algebra; cf. Theorem~\ref{thm:centertw}.
The matrix $S^*(u)=S^{-1}(-u-N/2)$ satisfies the quaternary relation \eqref{quatern}
and so, its Sklyanin determinant $\sdet S^*(u)$ is well-defined.

Fix a nonnegative integer $M\leq N$ such that $N-M$ is even and put $m=[M/2]$.
Set ${\cal A}=\{-n,\dots,-m-1,m+1\dots,n\}$
and ${\cal B}=\{-m,\dots,m\}$ and use the notation of Section~\ref{subsec:qdf}
for submatrices of $S(u)$.

\bth\label{thm:blocktw}
In the algebra $\wt{\Y}^{\pm}(N)$ we have
\ben
\sdet S(u)\ts\sdet S^*(-u+N/2-1)_{{\cal A}{\cal A}}= \sdet S(u)_{{\cal B}{\cal B}}.
\end{equation*}
\eth

We shall use the notation $\spr_{ij}(u)$ for the matrix elements of the
matrix $S^{-1}(u)$.

\bpr\label{prop:matelintw}
In the algebra $\wt{\Y}^{\pm}(N)$ 
the matrix elements of the matrices
$S(u)_{{\cal A}{\cal A}}$ and $S^{-1}(v)_{{\cal B}{\cal B}}$
commute with each other.
\epr

\subsection{The centralizer construction}\label{subsec:cctw}

Let $\g_n\subset\gl^{}_N$ be the classical 
Lie algebra of type $B_n$, $C_n$ or $D_n$,
as defined in Section~\ref{subsec:twdefrel}.
Fix an integer $m$ satisfying $0\leq m\leq n$ if $N=2n$,
and $-1\leq m\leq n$ if $N=2n+1$. Denote by 
$\g_m(n)$ the subalgebra of $\g_n$ spanned by
the elements $F_{ij}$ subject to the condition $m+1\leq |i|,|j|\leq n$.
Let $\Ar_m(n)$ denote the centralizer
of $\g_m(n)$ in the universal enveloping algebra $\Ar(n)=\U(\g_n)$.
In particular, $\Ar_0(n)$ (respectively, $\Ar_{-1}(n)$) is the center of $\Ar(n)$.
Let $\Ar(n)^0$ denote the centralizer of $F_{nn}$
in $\Ar(n)$ and let $\Ir(n)$ be the left ideal in $\Ar(n)$ generated by the
elements $F_{in}$, $i=-n,\dots,n$. Then $\Ir(n)^0=\Ir(n)\cap \Ar(n)^0$
is a two-sided ideal in $\Ar(n)^0$ and one has a vector space
decomposition
\ben
\Ar(n)^0=\Ir(n)^0\oplus \Ar(n-1).
\end{equation*}
Therefore the projection of $\Ar(n)^0$ onto $\Ar(n-1)$ 
with the kernel $\Ir(n)^0$ is an algebra homomorphism. Its
restriction to the subalgebra $\Ar_m(n)$ defines
a filtration preserving homomorphism
\ben
\pi_n: \Ar_m(n)\to \Ar_m(n-1)
\end{equation*}
so that one can define the algebra $\Ar_m$ as
the projective limit with respect to this sequence of homomorphisms
in the category of filtered algebras.

We denote by $\h_n$ the diagonal Cartan subalgebra
of $\g_n$, and by $\mathfrak{n}^+$ and $\mathfrak{n}^-$
the subalgebras spanned by the upper triangular and 
lower triangular matrices, respectively.
We identify $\U(\h_n)$ with the algebra of polynomial
functions on $\h_n^*$ and let $\lambda_i$ denote
the function which corresponds to $F_{ii}$.
For $i=1,\dots,n$ denote
\beql{rhocc}
\rho_{-i}=-\rho_i=\begin{cases}
i-1\quad&\text{for}\quad\g_n=\oa^{}_{2n},\\
i-\frac12\quad&\text{for}\quad\g_n=\oa^{}_{2n+1},\\
i\quad&\text{for}\quad\g_n=\spa^{}_{2n},
\end{cases}
\end{equation}
and set $l_i=-l_{-i}=\lambda_{i}+\rho_{i}$.	
We also set $\rho_0=1/2$ in the case
of $\g_{n}=\oa^{}_{2n+1}$.
Recall \cite[Section~7.4]{d:ae} that the image $\chi(z)\in\U(\h_n)$ 
of an element $z$ of the center of
$\Ar(n)$ under
Harish-Chandra isomorphism $\chi$ is uniquely determined
by the condition
\ben
z-\chi(z)\in \big(\mathfrak{n}^- \ts\Ar(n)+\Ar(n)\ts\mathfrak{n}^+\big).
\end{equation*}
If we identify $\U(\h_n)$ with the algebra of polynomials 
in	the variables $l_1,\dots,l_n$
then $\chi(z)$ belongs to the subalgebra ${\rm M}^*(n)$
of those polynomials $f=f(l_1,\dots,l_n)$ which are invariant
under the shifted action of the Weyl group.
More precisely, $f$ must be invariant
under all permutations of the variables and all
transformations $l_{i}\mapsto\pm l_{i}$, where in the
case of $\g_n=\oa_{2n}$ the number of `$-$' has to be even.
In the case of the minimum value of $m$ ($m=0$ or $m=-1$, respectively) 
the homomorphisms $\pi_n$ are interpreted as the specialization
homomorphisms $\pi_n: {\rm M}^*(n)\to {\rm M}^*(n-1)$ such that
\ben
\pi_n: f(\la_1,\dots,\la_n)\mapsto f(\la_1,\dots,\la_{n-1},0).
\end{equation*}
The corresponding projective limit in the category of filtered
algebras is an analog of the algebra of shifted symmetric functions
which is denoted by ${\rm M}^*$. The elements of ${\rm M}^*$ are well-defined
functions on the set of all sequences $\la=(\la_{1},\la_{2},\dots)$ which
contain only a finite number of nonzero terms. 
The following families of elements of ${\rm M}^*$ are analogs of power
sums, elementary symmetric functions, and complete symmetric functions:
\ben
\bal
p_m(\la)&=\sum_{k=1}^{\infty}(l_{k}^{\ts 2m}-\rho_{k}^{\ts 2m}),\qquad m=1,2,\ldots,\\
1+\sum_{m=1}^{\infty}e_m(\la)t^m&=\prod_{k=1}^{\infty}
\frac{1+l_{k}^{\ts 2}\ts t}{1+\rho_{k}^{\ts 2}\ts t},\\ 
1+\sum_{m=1}^{\infty}h_m(\la)t^m&=\prod_{k=1}^{\infty}
\frac{1-\rho_{k}^{\ts 2}\ts t}{1-l_{k}^{\ts 2}\ts t}.
\eal
\end{equation*}
Each of the families $\{p_m\}$,
$\{e_m\}$, $\{h_m\}$ can be taken as a system of algebraically independent
generators of the algebra ${\rm M}^*$.

\bpr\label{prop:shifttw}
The algebra $\Ar_0$ or $\Ar_{-1}$ in the case $N=2n$ or $N=2n+1$, respectively,
is isomorphic to the algebra of shifted
symmetric functions ${\rm M}^*$.
\epr

Consider the homomorphism \eqref{evaltw} and take its composition with the automorphism
of $\Y^{\pm}(N)$ given by Proposition~\ref{prop:sklauto} to yield another homomorphism
$\varphi_n:\Y^{\pm}(N)\to\U(\g_n)$. Set $M=2m$ or $M=2m+1$ depending on whether
$N=2n$ or $N=2n+1$. The image of the restriction of 
$\varphi_n$ to the subalgebra $\Y^{\pm}(M)$	is contained in the centralizer $\Ar_m(n)$.

\bth\label{thm:cctw}
The sequence of homomorphisms
$(\varphi_n\ |\ n\geq m)$ defines an algebra embedding
$
\varphi: \Y^{\pm}(M)\hookrightarrow\Ar_m.
$
Moreover, one has
an isomorphism
\ben
\Ar_m={\rm M}^*\ot \Y^{\pm}(M),
\end{equation*}
where $\Y^{\pm}(M)$ is identified with its image under the embedding $\varphi$.
\eth

\subsection{Commutative subalgebras}\label{subsec:cstw}

Fix an $N\times N$ matrix $C$ with entries in $\CC$
such that $C^{\tss t}=C$ or $C^{\tss t}=-C$, where the transposition is defined in
\eqref{transptw}.
Consider the
algebra $\Y^{\pm}(N)[[u^{-1}]]\ot (\End \CC^N)^{\ot N}$
and for any $1\leq k\leq N$ introduce its element
\ben
S(u,k)=\langle S_1,\ldots,S_k\rangle
\end{equation*}
as in Section~\ref{subsec:sklyanin} with the
variables $u_i$ specialized to $u_i=u-i+1$ for $i=1,\dots,k$.
Similarly, define the element $C(u,k)$ by
\ben
C(u,k)=C_{k+1}\wt R^{\tss t}_{k+1,k+2}\cdots \wt R^{\tss t}_{k+1,N}\ts
C_{k+2}\ts\wt R^{\tss t}_{k+2,k+3}\cdots \wt R^{\tss t}_{k+2,N}\cdots
C_{N-1}\ts\wt R^{\tss t}_{N-1,N}\ts C_N,
\end{equation*}
where we abbreviate $\wt R^{\tss t}_{ij}=R^{\tss t}_{ij}(-2u-N+i+j+2)$.
Introduce
the series $\sigma^{}_k(u,C)$ with coefficients in $\Y^{\pm}(N)$ by
\ben
\sigma^{}_k(u,C)=\tr A_N\ts S(u,k)\left(\ts 
\prod_{i=1,\dots,k}^{\rightarrow}\ts
\prod_{j=k+1,\dots,N}^{\rightarrow} R^{\tss t}_{ij} \right)
C(u,k),
\end{equation*}
where $R^{\tss t}_{ij}=R^{\tss t}_{ij}(-2u+i+j+2)$,
$A_N$ is the antisymmetrizer defined by \eqref{antisym} and the trace
is taken over all $N$ copies of $\End\CC^N$. 

\bth\label{thm:commsubtw}
All the coefficients of the series $\sigma^{}_1(u,C),\dots,\sigma^{}_N(u,C)$ commute
with each other.
Moreover, if the matrix $C$ has simple spectrum 
and satisfies $C^{\tss t}=-C$ then
these coefficients generate a 
maximal commutative subalgebra of the twisted Yangian $\Y^{\pm}(N)$.
\eth

Consider the epimorphism $\pi:\Y^{\pm}(N)\to \U(\g_n)$ defined in \eqref{evaltw}
(recall that $\g_n$ denotes the classical Lie algebra $\oa_{2n}$, $\oa_{2n+1}$
or $\spa_{2n}$).
Clearly, the coefficients of the images
of the series $\sigma^{}_k(u,C)$, $k=1,\dots,N$ under 
$\pi$ form a commutative subalgebra
${\mathcal C}\subseteq \U(\g_n)$. 

\bth\label{thm:commsububcd}
If the matrix $C$ has simple spectrum and satisfies 
$C^{\tss t}=-C$ then the subalgebra
${\mathcal C}$ of \ts $\U(\g_n)$ is maximal commutative.
\eth

\subsection*{Bibliographical notes}\label{subsec:bibbcd}

The twisted Yangians
were introduced by Olshanski in \cite{o:ty} where he also outlined 
their basic properties. 
A detailed exposition of the most of the results presented here
can be found in \cite{mno:yc}.

\noindent
{\bf \ref{subsec:sklyanin}.} Theorem~\ref{thm:sklproj} is proved in \cite{m:sd}.
Theorem~\ref{thm:stirling} was conjectured by Lascoux and proved by the author 
in \cite{m:sp}.

\noindent
{\bf \ref{subsec:sdf}.} Quasi-determinant factorization of the Sklyanin
determinant is given in \cite{m:ns}.

\noindent
{\bf \ref{subsec:cctw}.} The centralizer construction 
originates from Olshanski~\cite{o:ty}. A detailed proof
of Theorem~\ref{thm:cctw} is given in \cite{mo:cc}.

\noindent
{\bf \ref{subsec:cstw}.} The commutative subalgebras in the twisted Yangians
and the classical enveloping algebras
were constructed in Nazarov--Olshanski~\cite{no:bs}.

\medskip
Finite-dimensional irreducible representations of the twisted Yangians
are classified in \cite{m:fd} with the symplectic case done earlier in \cite{m:rt}.
Explicit constructions of all representations of $\Y^{\pm}(2)$ are also given in
\cite{m:fd}. These results were used in \cite{m:br,m:wb,m:wbg} to 
construct weight bases
of Gelfand--Tsetlin type for the classical Lie algebras.
Ragoucy~\cite{r:ty} discovered a relationship between the twisted Yangians and folded
${\mathcal W}$-algebras.
A family of algebras defined by a quaternary type relation
(or reflection equation) is defined by Sklyanin~\cite{s:bc}.
He also constructed commutative subalgebras and some representations of these
algebras. They were also studied in the physics literature in connection with
the integrable models with boundary conditions and the
nonlinear Schr\"odinger equation; see e.g. 
Kulish--Sklyanin~\cite{ks:as}, Kulish--Sasaki--Schwiebert~\cite{kss:cs},
Kuznetsov--J\o{rgensen}--Christiansen~\cite{kjc:nb},
Liguori--Mintchev--Zhao~\cite{lmz:be},
Mintchev--Ragoucy--Sorba--Zaugg~\cite{mrsz:ys,mrs:ss}. 
The action of such algebras
on hypergeometric functions was studied by Koornwinder and Kuznetsov~\cite{kk:gh}.

\section{Applications to classical Lie algebras}\label{sec:appl}
\setcounter{equation}{0}

Here we give constructions of families of Casimir elements for the
classical Lie algebras implied by the results discussed in the previous
sections. All of these constructions (including some well known)
are related with the quantum determinant for the Yangian $\Y(n)$
or the Sklyanin determinant for the twisted Yangian $\Y^{\pm}(N)$.
We keep the notation $\g_n$ for the classical Lie algebra
of type $B_n$, $C_n$ or $D_n$ as in Section~\ref{subsec:twdefrel}.
For any element $z$ of the center of the universal enveloping algebra
$\U(\gl_n)$ or $\U(\g_n)$ we shall denote by $\chi(z)$
its Harish-Chandra image; see Sections~\ref{subsec:cc} and \ref{subsec:cctw},
respectively. In the case of $\g_n$ we keep using the parameters $\rho_i$
defined in \eqref{rhocc}.

\subsection{Newton's formulas}\label{subsec:nf}

Consider the case of $\gl_n$ first. As in Section~\ref{subsec:mot}
we denote by $E$ the $n\times n$-matrix whose $ij$-th entry if $E_{ij}$.
Denote by $C(u)$ the {\it Capelli determinant\/}
\ben
C(u)=\sum_{p\in\Sym_n}\sgn p\cdot (u+E)_{p(1),1}\cdots (u+E-n+1)_{p(n),n}.
\end{equation*}
This is a polynomial in $u$ whose coefficients belong to the center
of $\U(\gl_n)$. The image of $C(u)$ under the Harish-Chandra isomorphism
is clearly given by
\beql{chicu}
\chi: C(u)\mapsto (u+l_1)\cdots (u+l_n),
\end{equation}
where $l_i=\la_i-i+1$.
Another family of central elements is provided
by the Gelfand invariants $\tr E^k$; see Section~\ref{subsec:mot}.

The following can be regarded as a noncommutative analog of the 
classical Newton formula which relates the elementary and power sums symmetric
functions; see e.g. Macdonald~\cite{m:sf}.

\bth\label{thm:newton}
We have the formula
\ben
1+\sum_{k=0}^{\infty}\frac{(-1)^k\ts \tr E^k}{(u-n+1)^{k+1}}=\frac{C(u+1)}{C(u)}.
\end{equation*}
\eth

\Proof Apply the homomorphism \eqref{eval} to the quantum Liouville formula; 
Section~\ref{subsec:qlf}. \endproof

\bco\label{cor:pp}
The images of the Gelfand invariants under the Harish-Chandra isomorphism
are given by
\beql{pp}
1+\sum_{k=0}^{\infty}\frac{(-1)^k\ts \chi(\tr E^k)}{(u-n+1)^{k+1}}=
\prod_{i=1}^n\big(1+\frac{1}{u+l_i}\big).
\end{equation}
\eco

In the case of $\g_n$ we denote by $F$ the $N\times N$-matrix whose
$ij$-th entry is $F_{ij}$. Introduce the {\it Capelli-type determinant\/}
\ben
C(u)=(-1)^n\sum_{p\in\Sym_N}
\sgn p\tss p'\cdot (u+\rho_{-n}+F)_{-a^{}_{p(1)},\ts a^{}_{p'(1)}}\cdots
(u+\rho_{n}+F)_{-a^{}_{p(N)},\ts a^{}_{p'(N)}},
\end{equation*}
where $(a_1,\dots,a_N)$ is any permutation of the indices $(-n,\dots,n)$ and
$p'$ is the image of $p$ under the map \eqref{mapsym}.
Using Theorem~\ref{thm:sklproj} and applying the homomorphism \eqref{evaltw}
we deduce that all coefficients of the polynomial $C(u)$ belong to the center of
$\U(\g_n)$.
If we take $(a_1,\dots,a_N)=(-n,\dots,n)$ then the image
of $C(u)$	under the Harish-Chandra isomorphism is easily found. For $N=2n$
we get
\ben
\chi: C(u)\mapsto \prod_{i=1}^n(u^2-l^2_i),
\end{equation*}
and for $N=2n+1$
\ben
\chi: C(u)\mapsto \big(u+\frac12\big)\prod_{i=1}^n(u^2-l^2_i).
\end{equation*}

\bth\label{thm:newtontw}
If $N=2n$ we have
\ben
1+\frac{2u+1}{2u+1\mp 1}\sum_{k=0}^{\infty}\frac{(-1)^k\ts \tr F^k}
{(u+\rho_n)^{k+1}}=\frac{C(u+1)}{C(u)},
\end{equation*}
where the upper sign is taken in the orthogonal case and the lower sign 
in the symplectic case.
If $N=2n+1$ then the same formula holds with
$C(u)$ replaced by
\ben
{\overline C}(u)=\frac{2u}{2u+1}\ts C(u).
\end{equation*}
\eth

\Proof Apply the homomorphism \eqref{evaltw} to the quantum Liouville formula; 
Section~\ref{subsec:liouvtw}. \endproof

\bco\label{cor:ppbcd}
The images of the Gelfand invariants under the Harish-Chandra isomorphism
are given by
\beql{ppbcd}
1+\frac{2u+1}{2u+1\mp 1}\sum_{k=0}^{\infty}\frac{(-1)^k\ts \chi(\tr F^k)}
{(u+\rho_n)^{k+1}}=
\prod_{i=-n}^n\big(1+\frac{1}{u+l_i}\big),
\end{equation}
where 
the zero index is skipped in the product if $N=2n$, while for $N=2n+1$
one should set $l_0=0$.
\eco

\subsection{Cayley--Hamilton theorem}\label{subsec:cht}

The polynomials $C(u)$ turn out to be the noncommutative
characteristic polynomials for the matrices $E$ and $F$.
Consider the case of $\gl_n$ first.

\bth\label{thm:chgln} 
We have the identities
\beql{chgln}
C(-E+n-1)=0\qquad\text{and}\qquad C(-E^{\tss t})=0.
\end{equation}
\eth

\Proof Applying the homomorphism \eqref{eval} to the relations
\eqref{comatrix} and \eqref{cotrans} and multiplying by the denominators
we get
\ben
C(u)=\wh C(u)\ts (u+E-n+1)\qquad\text{and}\qquad 
C(u)=\wh C^{\tss t}(u-1)\ts (u+E^{\tss t}),
\end{equation*}
where $\wh C(u)$ is polynomial in $u$ with coefficients in $\U(\gl_n)\ot\End\CC^n$.
\endproof

Taking the images of the identities \eqref{chgln} in a highest weight 
representation $L$ of $\gl_n$ with the highest weight $(\la_1,\dots,\la_n)$ we derive
the {\it characteristic identities\/}.

\bco\label{cor:char}
The image of the matrix $E$ in $L$ satisfies the identities
\ben
\prod_{i=1}^n (E-l_i-n+1)=0	\qquad\text{and}\qquad
\prod_{i=1}^n (E^{\tss t}-l_i)=0,
\end{equation*}
where $l_i=\la_i-i+1$.\endproof
\eco

Now turn to the case of the Lie algebras $\g_n$. As before, we consider the three cases
simultaneously.

\bth\label{thm:chbcd} 
We have the identity
\beql{chbcd}
C(-F-\rho_n)=0.
\end{equation}
\eth

\Outline Apply \eqref{evaltw} to the relation \eqref{sklco}
and multiply by the denominators. \endproof

Let $L$ be a highest weight representation of $\g_n$ with the highest weight
$(\la_1,\dots,\la_n)$ with respect to the basis elements $F_{11},\dots,F_{nn}$
of the Cartan subalgebra $\h_n$. The following are
the {\it characteristic identities\/} for $\g_n$ which are obtained
by taking the image of \eqref{chbcd} in $L$.

\bco\label{cor:charbcd}
The image of the matrix $F$ in $L$ satisfies the identities
\ben
\prod_{i=-n}^n (F-l_i+\rho_n)=0,
\end{equation*}
where $l_i=\la_i+\rho_i$. The zero index is skipped in the product if $N=2n$,
while for $N=2n+1$
one should set $l_0=\frac12$.
\endproof
\eco

\subsection{Graphical constructions of Casimir elements}\label{subsec:gc}

For $1\leq m\leq n$ denote by $E^{(m)}$ the
$m\times m$-matrix with the entries $E_{ij}$, where
$i,j=1,\dots,m$.
Let ${\mathcal E}^{(m)}$ denote the complete oriented graph with
the vertices $\{1,\dots,m\}$, the arrow from $i$ to $j$ is labelled by
the $ij$-th matrix element of the matrix $E^{(m)}-m+1$.
Then every path in the graph
defines a monomial in the matrix elements in a natural way.
A path from $i$ to $j$ is called {\it simple\/} if it does not
pass through the vertices $i$ and $j$ except for the beginning and the end of
the path.
Using this graph introduce the elements $\Lambda_k^{(m)}$, $S_k^{(m)}$,
$\Psi_k^{(m)}$ and $\Phi_k^{(m)}$ of
the universal enveloping algebra $\U(\gl_n)$ as follows.
For $k\geq 1$

$(-1)^{k-1}\Lambda_k^{(m)}$ is the sum of all monomials labelling
simple paths in ${\mathcal E}^{(m)}$ of length $k$ going from $m$ to $m$;

$S_k^{(m)}$ is the sum of all monomials labelling paths in
${\mathcal E}^{(m)}$ of length $k$ going from $m$ to $m$;

$\Psi_k^{(m)}$ is the sum of all monomials labelling paths in
${\mathcal E}^{(m)}$ of length $k$ going from $m$ to $m$, the coefficient of
each monomial being the length of the first return to $m$;

$\Phi_k^{(m)}$ is the sum of all monomials labelling paths in
${\mathcal E}^{(m)}$ of length $k$ going from $m$ to $m$, the coefficient of
each monomial being the ratio of $k$ to the number of returns to $m$.

\bth\label{thm:gcgln} The center of the algebra
$\U(\gl_n)$ is generated by the scalars and each of the
following families of elements
\ben
\bal
\Lambda_k=&\sum_{i_1+\cdots+i_n=k}\Lambda_{i_1}^{(1)}\cdots
\Lambda_{i_n}^{(n)},\\
S_k=&\sum_{i_1+\cdots+i_n=k}S_{i_1}^{(1)}\cdots
S_{i_n}^{(n)},\\
\Psi_k=&\sum_{m=1}^n\Psi_k^{(m)},\\
\Phi_k=&\sum_{m=1}^n\Phi_k^{(m)},
\eal
\end{equation*}
where $k=1,2,\dots,n$. Moreover, $\Psi_k=\Phi_k$ for any $k$, and 
the images of $\Lambda_k$, $S_k$ and
$\Psi_k$ under the Harish-Chandra isomorphism are, respectively,
the elementary, complete and power sums symmetric functions of degree $k$
in the variables $l_1,\dots,l_n$.
\eth

\Proof Consider the polynomial $C(u)$ introduced in Section~\ref{subsec:nf}
and set $\wt C(t)=t^n\ts C(t^{-1})$. Applying the homomorphism \eqref{eval}
to the decomposition \eqref{qdfactor} we obtain a decomposition
in the algebra of formal series with coefficients in $\U(\gl_n)$,
\ben
\wt C(t)=|1+tE^{(1)}|_{11}\cdots |1+t(E^{(n)}-n+1)|_{nn},
\end{equation*}
and the factors are permutable. By \cite[Proposition~7.20]{gkllrt:ns},
the elements introduced above are now interpreted as follows
\ben
1+\sum_{k=1}^{\infty}\Lambda_k^{(m)}t^k=|1+t(E^{(m)}-m+1)|_{mm},
\end{equation*}
\ben
1+\sum_{k=1}^{\infty}S_k^{(m)}t^k=|1-t(E^{(m)}-m+1)|_{mm}^{-1},
\end{equation*}
\ben
\sum_{k=1}^{\infty}\Psi_k^{(m)}t^{k-1}
=|1-t(E^{(m)}-m+1)|_{mm}\ts
\frac{d}{ dt}\ts|1-t(E^{(m)}-m+1)|_{mm}^{-1},
\end{equation*}
\ben
\sum_{k=1}^{\infty}\Phi_k^{(m)}t^{k-1}
=-\frac{d}{ dt}\log(|1-t(E^{(m)}-m+1)|_{mm}).
\end{equation*}
Due to the relations between the classical symmetric functions \cite{m:sf},
the second part of the theorem follows from \eqref{chicu}. \endproof

Similarly, in the case of $\g_n$ for any $1\leq m\leq n$ 
denote by $F^{(m)}$ the matrix with the entries $F_{ij}$, where
$i,j=-m,-m+1,\dots,m$ (the index $0$ is skipped if $N=2n$).					          
Consider the complete oriented graph ${\mathcal F}_m$ with
the vertices $\{-m,-m+1,\dots,m\}$, the arrow from $i$ to $j$ is labelled by
the $ij$-th matrix element of the matrix $F^{(m)}+\rho_m$.

Introduce the elements $\Lambda_{k}^{(m)}$, $\wt{\Lambda}_{k}^{(m)}$,
$S_{k}^{(m)}$, $\wt S_{k}^{(m)}$,
$\Phi_{k}^{(m)}$, $\wt{\Phi}_{k}^{(m)}$ of
the universal enveloping algebra $\U(\g_n)$ as follows: for $k\geq 1$

$(-1)^{k-1}\Lambda_{k}^{(m)}$ (resp. $-\wt{\Lambda}_{k}^{(m)}$)
is the sum of all monomials labelling
simple paths in ${\mathcal F}^{(m)}$
(resp. simple paths that do not pass through $-m$)
of length $k$ going from $m$ to $m$;

$S_{k}^{(m)}$ (resp. $(-1)^k\wt S_{k}^{(m)}$)
is the sum of all monomials labelling paths in
${\mathcal F}^{(m)}$ (resp. paths that do not pass through $-m$) 
of length $k$ going from $m$ to $m$;

$\Phi_{k}^{(m)}$ (resp. $(-1)^k\wt{\Phi}_{k}^{(m)}$)
is the sum of all monomials labelling paths in
${\mathcal F}^{(m)}$ (resp. paths that do not pass through $-m$)
of length $k$ going from $m$ to $m$, the coefficient of
each monomial being the ratio of $k$ to the number of returns to $m$.

\bth\label{thm:gcbcd}
Each of the
following families of elements is contained in the center of the algebra
$\U(\g_n)$:
\ben
\bal
\Lambda_{2k}=&\sum_{i_1+\cdots+i_{2n}=2k}\wt{\Lambda}_{i_1}^{(1)}
\Lambda_{i_2}^{(1)}\cdots \wt{\Lambda}_{i_{2n-1}}^{(n)}
\Lambda_{i_{2n}}^{(n)},\\
S_{2k}=&\sum_{i_1+\cdots+i_{2n}=2k}\wt{S}_{i_1}^{(1)}
S_{i_2}^{(1)}\cdots \wt{S}_{i_{2n-1}}^{(n)}
S_{i_{2n}}^{(n)},\\
\Phi_{2k}=&\sum_{m=1}^n(\wt{\Phi}_{2k}^{(m)}
+{\Phi}_{2k}^{(m)}),
\eal
\end{equation*}
where $k=1,2,\dots$\ts. Moreover, the images of
$(-1)^k\Lambda_{2k}$, $S_{2k}$ and
$\Phi_{2k}/2$ under the Harish-Chandra isomorphism are, respectively,
the elementary, complete and power sums symmetric functions of degree $k$
in the variables $l_1^2,\dots,l_n^2$.
\eth

\Outline The proof is the same as for Theorem~\ref{thm:gcgln}
with the use of Theorem~\ref{thm:sdetdec} and the homomorphism \eqref{evaltw}.
\endproof

\subsection{Pfaffians and Hafnians}\label{subsec:ph}

Suppose first that $\g_n$ is the orthogonal Lie algebra 
$\oa_{2n}$ or $\oa_{2n+1}$.	For any $1\leq k\leq n$ consider a subset
$I\subseteq \{-n,\dots,n\}$ of cardinality $2k$ so that the elements
of $I$ are $i_1<\cdots< i_{2k}$. The matrix $[F_{i_p,-i_q}]$
is skew-symmetric.
Introduce the corresponding
{\it Pfaffian} $\Pf F^I$ by
\ben
\Pf F^I=\frac{1}{2^k\tss k!}\sum_{\sigma\in\Sym_{2k}}
\sgn\sigma\cdot F_{i_{\sigma(1)},-i_{\sigma(2)}}\cdots 
F_{i_{\sigma(2k-1)},-i_{\sigma(2k)}}.
\end{equation*}
Given a subset $I$ set $I^*=\{-i_{2k},\dots,-i_1\}$ and denote
\ben
c_k=(-1)^k\sum_{|I|=2k} \Pf F^I\ts\Pf F^{I^*}, \qquad k\geq 1,
\end{equation*}
and $c_0=1$.
Consider the Capelli-type determinant $C(u)$ introduced in Section~\ref{subsec:nf}.

\bth\label{thm:pfdec}
The elements $c_k$ belong to the center of $\U(\g_n)$ 
and one has a decomposition
\ben
C(u)=\sum_{k=0}^n {c_k}\ts{(u^2-\rho_{1}^2)\cdots (u^2-\rho_{n-k}^2)}
\end{equation*}
if $N=2n$, and
\ben
C(u)=\big(u+\frac12\big)
\sum_{k=0}^n {c_k}\ts{(u^2-\rho_{1}^2)\cdots (u^2-\rho_{n-k}^2)}
\end{equation*}
if $N=2n+1$.
Moreover, the image of $c_k$ under the Harish-Chandra isomorphism is given by
\ben
\chi(c_k)=(-1)^k \sum_{i_1<\cdots < i_k} (l_{i_1}^2-\rho_{i_1}^2)\cdots 
(l_{i_k}^2-\rho_{i_k-k+1}^2).
\end{equation*}
\eth

Using the definition of the Sklyanin determinant one can write
a similar expression for the Capelli-type determinant in any realization of the
orthogonal or symplectic Lie algebra. Consider  
the realization of $\oa_{N}$ corresponding to the canonical symmetric form so that
the elements of $\oa_{N}$ are skew-symmetric matrices with respect to the usual
transposition. For this discussion only, we use more standard numbering
$1,\dots,N$ of the rows and columns
of such matrices. Here $F$ will denote the $N\times N$-matrix whose $ij$ entry is 
$F_{ij}=E_{ij}-E_{ji}$. The elements $c_k$ are now given by
\ben
c_k=\sum_{|I|=2k} (\Pf F^I)^2,
\end{equation*}
where $I=\{i_1,\dots,i_{2k}\}$ is a subset of $\{1,\dots,N\}$ and
$\Pf F^I$ is the Pfaffian of the skew-symmetric matrix
$[F_{i_p,i_q}]$. 
The Capelli-type determinant is given by
\ben
C(u)=
\sum_{p\in\Sym_{N}}\sgn p\tss p'\cdot (u+F+\sigma^{}_1)^{}_{p(1),\tss p'(1)}\cdots 
(u+F+\sigma^{}_N)^{}_{p(N),\tss p'(N)},
\end{equation*}
where $\sigma^{}_i=N/2-i$ for $i\leq n$ and $\sigma^{}_i=N/2-i+1$ for $i> n$. 
Theorem~\ref{thm:pfdec} holds in the same form with $\rho^{}_i$ replaced
by $\sigma^{}_{n-i+1}$ for every $i=1,\dots,n$.
Introduce also more standard determinant
\beql{stand}
D(u)=\sum_{p\in\Sym_{N}}
\sgn p\ts \cdot (u+F+m)_{p(1),\tss 1}(u+F+m-1)_{p(2),\tss 2}
\cdots (u+F-m+1)_{p(N),\tss N},
\end{equation}
where $m=N/2$.

\bth\label{thm:pfsqu}
The coefficients of the polynomial $D(u)$ are central in $\U(\oa_N)$
and the following decomposition holds
\beql{pfsqu}
D(u)=\sum_{k=0}^n {c_k}\ts{(u+m-k)(u+m-k-1)\cdots (u-m+k+1)}.
\end{equation}
\eth

In particular, we get the following
two analogs of the relation $(\Pf A)^2=\det A$ which holds for a numerical
skew-symmetric $2n\times 2n$-matrix $A$.

\bco\label{cor:pfdet} 
If $N=2n$ then 
\ben
(\Pf F)^2=C(0)=D(0).
\end{equation*}
\eco

Now consider the symplectic Lie algebra $\spa_{2n}$ and return to our
standard notation.
For any $k\geq 1$ consider a sequence $I$ of indices from $\{-n,\dots,n\}$
of cardinality $2k$ so that the elements
of $I$ are $i_1\leq\cdots\leq i_{2k}$. Set $\wt F_{ij}=\sgn i\cdot F_{ij}$.
Then we have $\wt F_{i,-j}=\wt F_{j,-i}$.
Introduce the corresponding
{\it Hafnian} $\Hf F^I$ by
\ben
\Hf F^I=\frac{1}{2^k\tss k!}\sum_{\sigma\in\Sym_{2k}}
\wt F_{i_{\sigma(1)},-i_{\sigma(2)}}\cdots 
\wt F_{i_{\sigma(2k-1)},-i_{\sigma(2k)}}.
\end{equation*}
For each $I$ let $f_{\pm 1},\dots,f_{\pm n}$ be the multiplicities of
the indices $\pm 1,\dots,\pm n$ in $I$, and let 
\ben
\sgn I=(-1)^{f_{-1}+\cdots+f_{-n}}.
\end{equation*}
Set $I^*=\{-i_{2k},\dots,-i_1\}$ and denote
\ben
d_k=(-1)^k\sum_{|I|=2k} \frac{\sgn I\cdot\Hf F^I\ts\Hf F^{I^*}}
{f_1!f_{-1}!\cdots f_n!f_{-n}!}.
\end{equation*}
Consider the Capelli-type determinant $C(u)$ introduced in Section~\ref{subsec:nf}
and introduce the series $c(u)$ by
\ben
c(u)=\frac{C(u)}{(u^2-\rho^2_{1})\cdots (u^2-\rho^2_{n})}.
\end{equation*}

\bth\label{thm:hfdec}
The elements $d_k$ belong to the center of $\U(\spa_{2n})$ 
and one has a decomposition
\ben
c(u)^{-1}=1+\sum_{k=1}^{\infty} \frac{d_k}{(u^2-(n+1)^2)\cdots (u^2-(n+k)^2)}.
\end{equation*}
Moreover, the image of $d_k$ under the Harish-Chandra isomorphism is given by
\ben
\chi(d_k)=\sum_{i_1\leq\cdots \leq i_k} (l_{i_1}^2-i_1^2)\cdots 
(l_{i_k}^2-(i_k+k-1)^2).
\end{equation*}
\eth

\subsection*{Bibliographical notes}\label{subsec:appl}

{\bf \ref{subsec:nf}.} The images of the Gelfand invariants
under the Harish-Chandra isomorphism were first found by
Perelomov and Popov~\cite{pp:cou, pp:coos, pp:co}; see also
\v{Z}elobenko~\cite{z:cl}.
The observation that the formulas \eqref{pp} and \eqref{ppbcd}
are related with the theory of Yangians is due to Cherednik~\cite{c:qg}.
A derivation of the Perelomov--Popov formulas from the Liouville formula
is contained in \cite{m:yl}. Of course, given the Harish-Chandra images of the
polynomials $C(u)$, Theorems~\ref{thm:newton} and \ref{thm:newtontw}
follow from Corollaries~\ref{cor:pp} and \ref{cor:ppbcd}. Different proofs
of Newton's formulas without using the Yangians are given by Itoh
and Umeda~\cite{i:en, i:ce, u:nf} for the general linear and orthogonal
Lie algebras.

\noindent{\bf \ref{subsec:cht}.} 
Theorem~\ref{thm:chgln} is due to 
Nazarov and Tarasov~\cite{nt:yg}. Theorem~\ref{thm:chbcd} is proved in
\cite{m:sd}. An independent proof is given by Nazarov (unpublished)
without using the map \eqref{mapsym}. One more proof in the
orthogonal case is given by Itoh~\cite{i:ce} with the use of a determinant
of type \eqref{stand}
instead of $C(u)$.
Corollaries~\ref{cor:char} and \ref{cor:charbcd} are the remarkable
characteristic identities which are due to Bracken and Green~\cite{bg:vo, gr:ci}.
More general identities are obtained by Gould~\cite{g:ci}.

\noindent{\bf \ref{subsec:gc}.}
The elements $\La_k$, $S_k$, $\Psi_k$, $\Phi_k$ are the noncommutative
symmetric functions associated with a matrix; see 
Gelfand {\it et al\/}~\cite{gkllrt:ns}.
Theorem~\ref{thm:gcgln} is proved in \cite[Section~7.5]{gkllrt:ns}.
Theorem~\ref{thm:gcbcd}	is contained in \cite{m:ns}. A different version
for the orthogonal Lie algebra is given in \cite{m:sp}.

\noindent{\bf \ref{subsec:ph}.}
Most of these results are contained in \cite{mn:ci}; see also \cite{m:sp}.
The centrality of the determinant \eqref{stand} was first proved by
Howe and Umeda~\cite{hu:ci}. Its relationship with the Pfaffians
(formula \eqref{pfsqu})
was established by Itoh and Umeda~\cite{iu:ce}. The 
first relation in Corollary~\ref{cor:pfdet} 
is proved in \cite{m:sd}.
Both Pfaffians and Hafnians are key ingredients
in the analogs of the celebrated Capelli identity \cite{c:zc, c:ot}
for the classical Lie algebras given in \cite{mn:ci};
see also Weyl~\cite{w:cg}, Howe~\cite{h:rc}, Howe--Umeda~\cite{hu:ci} 
for the role of the Capelli identity in the classical invariant theory.
The polynomials $\chi(c_k)$ and $\chi(d_k)$ are respectively 
the elementary and complete factorial symmetric functions;
see e.g. Okounkov--Olshanski~\cite{oo:ss, oo:ss2}.

\end{document}